\documentclass[a4paper,11pt]{article}
\usepackage{amsmath,amscd,amssymb}
\usepackage{graphicx,psfrag}
\usepackage{xypic}

\def \vs {\vskip}
\def \hs {\hskip}
\def \noi {\noindent}

\def \oo {{\cal O}}

\def \a {{\alpha}}
\def \b {{\beta}}
\def \ga {{\gamma}}
\def \Ga {{\Gamma}}
\def \vp {\varphi}

\def \p {{\mathbb P}}
\def \Z {{\mathbb{Z}}}
\def \N {\mathbb{N}}

\def \cC {{\mathfrak{C}}}

\def \fl {\rightarrow}

\def \im {{\rm Im}}
\def \pic {{\rm Pic}}
\def \pu {{\mathbb P}^1}
\def \opu {{\oo_\pu}}
\def \Mor #1#2{{\bf{Hom}}_{#1}(\pu,#2)}

\def \dm {{\textsc{Proof. }}}

\def \tha #1#2{\noi{\bf#1{\uppercase{\footnotesize{#2}}}}}

\newtheorem{theor}{\tha{T}{heorem}}[section]
\newenvironment{theo}{
  \begin{theor}\hs -0.2 cm {\bf .} ---  }
{  \end{theor}}

\newtheorem{propo}[theor]{\tha{P}{roposition}}
\newenvironment{prop}{
  \begin{propo}\hs -0.2 cm {\bf .} ---  }
{  \end{propo}}

\newtheorem{lemma}[theor]{\tha{L}{emma}}
\newenvironment{lemm}{
  \begin{lemma}\hs -0.2 cm {\bf .} ---  }
{  \end{lemma}}

\newtheorem{fait}[theor]{\tha{F}{act}}
\newenvironment{fact}{
  \begin{fait}\hs -0.2 cm {\bf .} ---  }
{  \end{fait}}

\newtheorem{defini}[theor]{\tha{D}{efinition}}
\newenvironment{defi}{
  \begin{defini}\hs -0.2 cm {\bf .} ---  }
{  \end{defini}}

\newtheorem{corollaire}[theor]{\tha{C}{orollary}}
\newenvironment{coro}{
  \begin{corollaire}\hs -0.2 cm {\bf .} ---  }
{  \end{corollaire}}

\newtheorem{exemple}[theor]{\sc{Example}}
\newenvironment{exem}{
  \begin{exemple}\hs -0.2 cm {\bf .} ---  }
{  \end{exemple}}

\newtheorem{remarq}[theor]{\sc{Remark}}
\newenvironment{rema}{
  \begin{remarq}\hs -0.2 cm {\bf .} ---  }
{  \end{remarq}}

\newtheorem{preuve}{\sc{Proof}}

\def \comp {{\mathfrak{ne}}}
\def \l {{\lambda}}
\def \sca #1#2{\left\langle#1,#2\right\rangle}
\def \scal #1#2{\langle #1,#2 \rangle}
\def \Bt {{\widetilde{B}}}
\def \Xt {{\widetilde{X}}}
\def \ft {{\widetilde{f}}}
\def \Ct {{\widetilde{C}}}
\def \at {{\widetilde{\a}}}
\def \bt {{\widetilde{\b}}}
\def \Ch {{\widehat{C}}}
\def \fc {{\widehat{f}}}

\def \pd {\varpi}
\def \Del {\Delta}

\begin{document}

~
\vs -0.5 cm

\centerline{\large{\uppercase{\bf{Rational curves on }}}}
\centerline{\large{\uppercase{\bf{minuscule Schubert varieties}}}}

\vs 1 cm

\centerline{\Large{Nicolas \textsc{Perrin}}}

\vs 1.5 cm

\centerline{\Large{\bf Introduction}} 

\vs 0.5 cm

Let us denote by $\cC$ the variety of lines in $\p^3$ meeting a fixed
line, it is a grassmannian (and hence minuscule) Schubert variety. In
\cite{PE3} we described the irreducible components of the scheme of
morphisms from $\pu$ to $\cC$ and the general morphism of these
irreducible components.

\vs 0.2 cm

In this text we study the scheme of morphisms from $\pu$ to any
minuscule Schubert variety $X$. Let us recall that we studied in
\cite{PE} the scheme of morphisms from $\pu$ to any homogeous
variety. The main idea, in the case of a minuscule Schubert variety
$X$, is to restrict ourselves to the dense orbit under the
stabilisator ${\rm Stab(X)}$ of $X$ and apply the results of
\cite{PE}.

\vs 0.2 cm

More precisely, let $U$ be the dense orbit under ${\rm Stab(X)}$ in
$X$ and let $Y$ be the complementary. Because $X$ is a minuscule
Schubert variety the closed subset $Y$ of $X$ is of codimension at
least 2 (see paragraph \ref{divschub}). This fact and the
stratification of $X$ by Schubert subvarieties gives us a surjective
morphism (see paragraph \ref{preliminaires}):
$$s:\pic(U)^\vee\to A_1(X).$$
For any class $\a\in A_1(X)$, we can consider the following morphism:
$$i:\!\!\!\coprod_{s(\b)=\a}\Mor{\b}{U}\fl\Mor{\a}{X}$$ 
where $\Mor{\a}{X}$ is the scheme of morhisms $f:\pu\to X$ with
$f_*[\pu]=\a$ and $\Mor{\b}{U}$ is the scheme of morhisms $g:\pu\to U$
such that $[g]=\b$ where $[g]$ is the linear function $L\mapsto{\rm
  deg}(g^*L)$ on $\pic(U)$. As $Y=X\setminus U$ lies in codimension 2,
we expect the image of this morphism to be dense (this is the crucial
point of the proof). This condition means that any morphism $\pu\to X$
can be deformed such that the image of this deformation does not meet
$Y$. If the morphism $i$ defined above is dominant, we may apply the
results of \cite{PE} to prove that $\Mor{\b}{U}$ is irreducible as
soon as it is non empty and the images of these irreducible
$\Mor{\b}{U}$ will give the irreducible components of $\Mor{\a}{X}$.

Let us denote by $\comp(\a)$ the subset of $\pic(U)^\vee$ given elements
$\b$ such that $s(\b)=\a$ and $\Mor{\b}{U}$ is non empty (see
paragraph \ref{preliminaires} for a more precise definition in terms
of roots). We prove the

\vs 0.4 cm

\begin{theo}
The irreducible components of the scheme of
  morphisms $\Mor{\a}{X}$ are indexed by $\mathfrak{ne}(\a)$.
\end{theo}

Here is an outline of the paper. In the first paragraph we define the
surjective map $s$ of the introduction and the set $\comp(\a)$ for $X$
any Schubert variety and $\a\in A_1(X)$. In the second paragraph we
recall the definition of a minuscule Schubert variety and its
properties. We also prove a positivity result on roots we will need
later. In the third paragraph we recall the construction of the
Bott-Samelson resolution $\pi:\Xt\to X$ of a Schubert variety $X$ and
describe some cycles on $\Xt$. In the fourth paragraph, we construct
some big families of  curves on $\Xt$ contracted by $\pi$. In the
fifth paragraph we study the scheme of morphisms $\Mor{\at}{\Xt}$ and
prove some smoothing results with the curves contructed in the fourth
paragraph. In the last paragraph we prove our main result.

The key point as indicated above is to prove that the map $i$ is
dominant that is to say that any morphism $f:\pu\to X$ can be
factorised in $U$ (modulo deformation). We prove this by lifting $f$
in $\ft$ on $\Xt$. It is 
now sufficent to prove that the lifted curve $\ft$ of a general curve
$f$ does not meet the divisors contracted by $\pi$. If $\ft$ does meet
a contracted divisor $D$ then we add a "line" $L\subset D$ with
$L\cdot D=-1$ constructed in the fourth paragraph and smooth the union
$\ft(\pu)\cup L$. The intersection with $D$ is lowered by  one in the
operation. We conclude by induction on the number of
intersection of $\ft$ with the contracted divisors.

\begin{rema}
  {\rm (\i)} The variety $\cC$ can also be seen as a cone over a smooth
2-dimensional quadric embedded in $\p^3$. We treat more generaly the
case of a cone $X$ over an homogeneous variety in the forthcoming
paper \cite{PE2}. In this situation we can also define for $\a\in
A_1(X)$ a class $\comp(\a)$ as previously but the irreducible
components of $\Mor{\a}{X}$ are not always indexed by
$\comp(\a)$. It is the case if and only if the projectivised tangent
cone of the singularity (here the embbeded homogeneous variety)
contains lines. 

{\rm (\i\i)} This condition on the existence of lines in the
projectivised tangent cone of the singularity also appears for more
general Schubert varieties.

{\rm (\i\i\i)} In \cite{Brionpolo}, M. Brion and P. Polo proved that
the singularities of minuscule Schubert varieties are locally
isomorphic to cones over homogeneous varieties. With the results of
\cite{PE2} this implies that the key problem of factorising morphisms
trough $U$ is locally true. Unfortunatly it is not obvious to prove
the global results thanks to this local property. It is nevertheless a
good guide for intuition and we solve here the global problem using
Bott-Samelson resolutions. 

\end{rema}

\section{Preliminary}
\label{preliminaires}

In this paragraph we explain the results on cycles used in the
introduction. We describe the surjective morphism $s:\pic(U)^\vee\to
A_1(X)$ and define the set of classes $\comp(\a)$ for $\a\in A_1(X)$.

\vs 0.2 cm

Let $X$ be a scheme of dimension $n$. 
Denote by $Z_*(X)$ the group of 1-cycles on  $X$ and by
$Z^\equiv_*(X)$ %, $Z_*^{a}(X)$ 
and $Z^{r}_*(X)$ the subgroups of cycles trivial for the numerical
%, algebraic 
and rational equivalence. Let us denote by $N_*(X)$ %, $B_*(X)$ 
and $A_*(X)$ the corresponding quotients. The Picard group is the
image in $A_{n-1}(X)$ of the subgroup of Cartier
divisors in $Z_{n-1}(X)$ and we denote by $N^1(X)$ the quotient of
$\pic(X)$ by numerical equivalence.

\begin{lemm}
  Let $X\subset G/P$ be a Schubert variety ($G$ a Lie group and $P$ a
  parabolic subgroup of $G$). Then one has

{\rm (\i)} $\pic(X)\simeq N^1(X)$,

{\rm (\i\i)}  $A_1(X)\simeq N_1(X)$. 

\noi
In particular we have $A_1(X)\simeq\pic(X)^\vee$. 
\end{lemm}

\dm
(\i) Thanks to the results of \cite{FS..} the groups $A_*(X)$ are free
generated by Schubert subvarieties and furthermore rational and
algebraic equivalence are the same. So on the one hand, the Picard
group is contained in $A_{n-1}(X)$ and is in particular free. 

On the other hand, thanks to \cite{F} Example 19.3.3, we know that a
Cartier divisor $D$ is numerically trivial if for some $m\in\N$ we
have $mD$ is algebraically trivial. This implies for Schubert
varieties that $mD$ is rationaly trivial and because $\pic(X)$ is
torsion free $D$ is trivial in $\pic(X)$. This implies that
$\pic(X)\simeq N^1(X)$.

(\i\i) The results of \cite{FS..} also imply that $A_1(X)$ is
generated by the one-dimensional Schubert varieties in
$X$. But on $G/P$ there is a duality between the Picard group and one
dimensional Schubert varieties. In praticular for any one
dimensional Schubert variety $Z$ there is a line bundle $L_Z$ such that
$L_Z\cdot Z=1$ and $L_Z$ is trivial on any other  one
dimensional Schubert variety. If the $Z_i\subset X$ are the one
dimensional Schubert varieties in $X$ then the restrictions of
the $L_{Z_i}$ to $X$ form a dual family to the $Z_i$. In particular
the $Z_i$ are numerically independent. As they form a
basis of $A_1(X)$ we have $A_1(X)\simeq N_1(X)$.

%The Picard group is contained in $A_{n-1}(X)$ and is in particular
%free. Thanks to \cite{F} this implies thatfree generated by a very
%ample divisor. It 
%is clearly dual to $N_1(X)$ thanks to the degree morphism. 

The duality
comes from general duality between $N_1(X)$ and $N^1(X)$.\hfill$\Box$

\vs 0.4 cm

Let $U$ be the smooth locus of $X$. If $X$ is minuscule (see
definition in paragraph \ref{minuscule}) this smooth locus $U$ is the
dense orbit under ${\rm Stab(X)}$ in $X$ (see
\cite{Brionpolo}\footnote{We do not need the results of \cite{Brionpolo}
  to define $\comp(\a)$, see theorem \ref{final}, but it is more simple
  with this fact on the singular locus.}). Let $Y$ be the complementary
of $U$ in $X$. Because $X$ is a normal variety the closed subset $Y$
is of codimension at least 2, this in particular implies that
$\pic(U)=A_{n-1}(U)\simeq A_{n-1}(X)$. We now have the following
inclusion: 
$$\pic(X)\subset A_{n-1}(X)\simeq\pic(U)$$
giving the surjection
$$s:\pic(U)^\vee\to A_1(X).$$
With these notations we make the following:

\begin{defi}
  Let $X$ be any Schubert variety and let $\a\in A_1(X)$. We define
  the set $\comp(\a)\subset A_{n-1}(X)^\vee$. 

Let us make the
  identification $A_{n-1}(X)\simeq\pic(U)$. The elements of
  $\comp(\a)$ are the elements $\b\in\pic(U)^\vee$ such that
  $s(\b)=\a$ and there exists a curve $C\subset U$ with $[C]=\b$
  as a linear form on $\pic(U)$ ($\b$ is effective). 
\end{defi}

In the case of minuscule Schubert variety $X\subset G/P$ we describe
$\comp(\a)$ more precisely: the smooth part $U$ is the dense orbit
under ${\rm Stab X}$. Let $R$ be the levi subgroup of ${\rm
  Stab}(X)$, the orbit $U$ is of the form $QP/P\simeq Q/Q\cap P$ where
$Q={\rm Stab}(X)$ is a 
parabolic subgroup of $G$. We proved in \cite{PE} proposition 5 that
this orbit is a tower of affine bundles over the homogeneous variety
$R/R\cap P$. In particular $\pic(U)\simeq\pic(R/R\cap P)$ is given in terms
of weights with a particular weight given by the generator of
$\pic(X)$. Furthermore we proved in \cite{PE} that the elements
$\b\in\pic(R/R\cap P)^\vee$ are effective if they are in the dual cone
of the cone of effective divisor, in other words they correpond to
positive roots. 

\begin{exem}
  If $X$ is a grassmannian Schubert variety given by a partition
  $\l$, consider the associated Young diagram (see for example
  \cite{maniv}). Then the Picard group $\pic(U)$ is free
  and has as many generators $(L_i)_{i\in[1,r]}$ as the numbers of
  holes. The generator $L$ of $\pic(X)$ is given by 
$$L=\sum_{i\in[1,r]}L_i.$$
If $\a\in A_1(X)$ is such that $\a\cdot L=d$ then $\comp(\a)$ is given
by the $r$-tuples $(b_i)_{i\in[1,r]}$ of non negative integers such that 
$$\sum_{i\in[1,r]}b_i=d.$$
The number of irreducible components is $\binom{d+r-1}{d}$.
\end{exem}

\begin{rema}
  The scheme $\Mor{\a}{X}$ is the scheme of
morphisms from $\pu$ to $X$ of class $\a$ (for more details see
\cite{GR} and \cite{MO}). 

In general, this will just mean that $\a\in A_1(X)$
  and that $f_*[\pu]=\a$ but sometimes (in particular in the
  introduction for the open part $U$) we consider $\a\in\pic(X)^\vee$
  and the class of a morphism $f:\pu\to X$ will be the linear form
  $\pic(X)\to\Z$ given by $L\mapsto{\rm deg}(f^*L)$.

In the case of a minuscule Schubert variety $X$ the two notion
coincide because of the previous lemma. 

In the case of the open part
$U$ of a minuscule Schubert variety $X$, these scheme are connected
components of the scheme of morphisms with a fixed 1-cycle class.
\end{rema}

\section{Minuscule Schubert varieties}
\label{minuscule}

\subsection{Definitions}

In this paragraph we recall the notion
of minuscule weight and study the related homogeneous
and Schubert varieties. Our basic reference will be \cite{GofG/P3}.

Let $G$ be a semi-simple algebraic group, fix $T$ a maximal torus and
$B$ a Borel subgroup containing $T$. Let us denote by $\Del$ the set
of all roots, by $\Del^+$ (resp. $\Del^-$) the set of positive
(resp. negative) roots, by $S$ the set of simple roots associated to
the data $(G,T,B)$ and by $W$ the associated Weyl group. If $P$ is a
parabolic subgroup containing $B$ we note $W_P$ the subgroup of $W$
corresponding to $P$. Let us finally denote by $\Bt$ the opposite Borel
subgroup (corresponding to the negative roots) and by $i$ the Weyl
involution on simple roots. This involution sends a simple root $\b$ on
$-w_0(\b)$ and is also defined on fundamental weights.

\begin{defi}
Let $\pd$ be a fundamental weight, 

(\i) we say that $\pd$ is minuscule if we have
$\sca{\a^\vee}{\pd}\leq1$ for all positive root $\a\in\Del^+$;

(\i\i) we say that $\pd$ is cominuscule if $\sca{\a_0^\vee}{\pd}=1$
where $\a_0$ is the longuest root.
\end{defi}

With the notation of N. Bourbaki \cite{bourb}, the minuscule and
cominuscule weights are:

%\begin{center}
%\begin{tabular}{|c|c|c|}
%\hline
%Type&minuscule\\
%\hline
%$A_n$&$\pd_1\cdots\pd_n$\\
%\hline
%$B_n$&$\pd_n$\\
%\hline
%$C_n$&$\pd_1$\\
%\hline
%$D_n$&$\pd_1$, $\pd_{n-1}$ and $\pd_n$\\
%\hline
%$E_6$&$\pd_1$ and $\pd_6$\\
%\hline
%$E_7$&$\pd_7$\\
%\hline
%$E_8$&none\\
%\hline
%$F_4$&none\\
%\hline
%$G_2$&none\\
%\hline
%\end{tabular}
%\end{center}

\begin{center}
\begin{tabular}{|c|c|c|}
\hline
Type&minuscule&cominuscule\\
\hline
$A_n$&$\pd_1\cdots\pd_n$&same weights\\
\hline
$B_n$&$\pd_n$&$\pd_1$\\
\hline
$C_n$&$\pd_1$&$\pd_n$\\
\hline
$D_n$&$\pd_1$, $\pd_{n-1}$ and $\pd_n$&same weights\\
\hline
$E_6$&$\pd_1$ and $\pd_6$&same weights\\
\hline
$E_7$&$\pd_7$&same weight\\
\hline
$E_8$&none&none\\
\hline
$F_4$&none&none\\
\hline
$G_2$&none&none\\
\hline
\end{tabular}
\end{center}

\vs 0.2 cm

\begin{rema}
  The Weyl involution $i$ acts on minuscule and on cominuscule weights.
\end{rema}

\begin{defi}
  Let $\pd$ be a  minuscule weight and let $P_\pd$
be the associated parabolic subgroup. The homogeneous variety
$G/P_\pd$ is then said to be minuscule. The Schubert varieties of a
minuscule homogeneous variety are called minuscule Schubert
varieties.
\end{defi}

\begin{rema}
  To study minuscule homogeneous varieties and their
  Schubert varieties, it is sufficent to restrict ourselves to
  simply-laced groups.

In fact the variety $G/P_{\pd_n}$ with $G={\rm Spin}_{2n+1}$ is
isomorphic to the variety $G'/P'_{\pd_{n+1}}$ with
$G'={\rm Spin}_{2n+2}$ and there is a one to one correspondence
between Schubert varieties thanks to this isomorphism. The same
situation occurs with $G/P_{\pd_1}$, $G={\rm Sp}_{2n}$ and
$G'/P'_{\pd_1}$, $G'={\rm SL}_{2n}$.
\end{rema}

\subsection{Divisors on minuscule Schubert varieties}
\label{divschub}

In this paragraph we describe the divisors on minuscule Schubert
varieties. For proofs and more details see \cite{GofG/P3}. 

\begin{defi}
\label{moving}
  Let $\bar\phi\in W/W_{P_{\pd}}$ and let $X(\bar\phi)$ the associated
  Schubert variety. A Schubert divisor $X(s_\b\bar\phi)$ in
  $X(\bar\phi)$ defined by a \emph{simple} root $\b$ is called a
  moving divisor. All other Schubert divisor are said to be
  stationary.
\end{defi}

\begin{rema}
  The term "moving divisor" comes from the fact that the
  Schubert variety $X(\bar\phi)$ is stable under the action of
  $U_{-\b}$ whereas $X(s_\b\bar\phi)$ is moved by $U_{-\b}$ in
  $X(\bar\phi)$ (see \cite{GofG/P3}).
\end{rema}

We have the following proposition (\cite{LW} Lemma 1.14):

\begin{prop}
  With the notation of the definition \ref{moving} then
  $X(s_\b\bar\phi)$ is a moving divisor in $X(\bar\phi)$ if and only
  if $\bar\phi$ has a reduced expression starting with $s_\b$.
\label{movingdeb}
\end{prop}

We now have the following theorem (\cite{kempf}, Th. 1 or
\cite{GofG/P3} Th. 3.10) which describes the divisors of a minuscule
Schubert variety:

\begin{theo} 
  Let $X$ be a minuscule Schubert variety,
  then every Schubert divisor in $X$ is a moving divisor.
\label{movingmin}
\end{theo} 

\begin{rema}
  {\rm (\i)} This theorem is equivalent to the fact that weak and
  strong Bruhat orders coincide on minuscule Schubert varieties.

{\rm (\i\i)} Let $U$ be the dense orbit in $X$ under the action of
stabilisator ${\rm Stab}(X)\subset G$. Let $Y$ be the complementary of
$U$ in $X$. A consequence of this theorem is that $Y$ is in
codimension at least 2.
\end{rema}

\subsection{A positivity result}

Let $(\ga_i)_{i\in[1,n]}$ be a sequence of simple roots and define
$\phi=s_{\ga_1}\cdots s_{\ga_n}$. We suppose in addition that
$l(\phi)=n$. Set $\b_i=i(\ga_i)$ and let us define a sequence of roots
$(\a_i)_{i\in[1,n]}$ by
$$\a_1=\b_1,\ \a_2=s_{\b_1}(\b_2),\ \dots,\
\a_n=s_{\b_1}\cdots s_{\b_{n-1}}(\b_n).$$
Remark that this construction is symetric in the sense that if the
$(\a_i)_{i\in[1,n]}$ are given we can recover the $(\b_i)_{i\in[1,n]}$
by the formulae
$$\b_1=\a_1,\ \b_2=s_{\a_1}(\a_2),\ \dots,\
\b_n=s_{\a_1}\cdots s_{\a_{n-1}}(\a_n).$$

\begin{rema}
  We use these notations to fit with those of the
  Bott-Samelson resolution.
\end{rema}

\begin{prop}
  \label{posi}
  Let $\pd$ be a minuscule weight. Suppose that $\phi$ is the smallest
  element in the class $\bar\phi\in W/W_{P_\pd}$. Then for all
  $(i,j)\in[1,n]$ we have
  $$\sca{\a_i^\vee}{\a_j}\geq0.$$
\end{prop}

\dm
Let us define the sequence $(\bt_i)_{i\in[1,n]}$ of simple roots as
beeing the sequence $(\b_i)_{i\in[1,n]}$ with reversed order, that is
to say $\bt_i=\b_{n+1-i}$. With this sequence we can construct a
sequence $(\at_i)_{i\in[1,n]}$ by
$$\at_1=\bt_1,\ \at_2=s_{\bt_1}(\bt_2),\ \dots,\
\at_n=s_{\bt_1}\cdots s_{\bt_{n-1}}(\bt_n).$$

\begin{lemm}
  For all $i\in[1,n]$, we have
$$\sca{\a_i^\vee}{\a_j}=\sca{\at_{n+1-i}^\vee}{\at_{n+1-j}}.$$
\end{lemm}

\dm
we have
$$\sca{\at_{n+1-i}^\vee}{\at_{n+1-j}}=\sca{s_{\bt_1}\cdots
  s_{\bt_{n-i}}(\bt_{n+1-i})^\vee}{s_{\bt_1}\cdots
  s_{\bt_{n-j}}(\bt_{n+1-j})}\ \ \ \ \ \ \ \ \ \ \ \ \ \ \ \ \ \
\ \ \ \ \ \ \ \ \ \ \ \ \ \ \ \ \ \ \ \ \ \ \ \ \ \ \ \ \ \ \ \ \ \ \ $$
$$\ \ \ \ \ \ \ \ \ \ \ \ \ \ \ \ =\sca{s_{\b_n}\cdots
  s_{\b_{i+1}}(\b_{i})^\vee}{s_{\b_n}\cdots
  s_{\b_{j+1}}(\b_{j})}=\sca{s_{\b_1}\cdots
  s_{\b_{i}}(\b_{i})^\vee}{s_{\b_1}\cdots
  s_{\b_{j}}(\b_{j})}$$
where we applied $s_{\b_1}\cdots s_{\b_{n}}$ to get the last
  equality. But we have
$$\sca{s_{\b_1}\cdots
  s_{\b_{i}}(\b_{i})^\vee}{s_{\b_1}\cdots
  s_{\b_{j}}(\b_{j})}=\sca{s_{\b_1}\cdots
  s_{\b_{i-1}}(-\b_{i})^\vee}{s_{\b_1}\cdots
  s_{\b_{j-1}}(-\b_{j})}=\sca{\a_i^\vee}{\a_j}.$$
\vs -0.9 cm \hfill$\Box$

\vs 0.4 cm

It is thus enough to prove the result on the sequence
$(\at_i)_{i\in[1,n]}$. As $\phi$ is the smallest element in $\bar
\phi$, the reduced expression $\phi=s_{\ga_1}\cdots
s_{\ga_n}=s_{i(\bt_n)}\cdots s_{i(\bt_1)}$ in $W$ is still
reduced in $W/W_{P_\pd}$. Let us prove the following lemma:

\begin{lemm}
\label{pgqb}
Let $\b$ the only simple root such that $\sca{\b^\vee}{i(\pd)}=1$. For
all $i\in[1,n]$, the roots $\at_i$ are such that
$$\at_i\geq\b.$$
\end{lemm}

\dm
Because the expression $\bar\phi=s_{i(\bt_n)}\cdots s_{i(\bt_1)}$
  is reduced in $W/W_{P_\pd}$, we have for all $i\in[1,n]$
$$\sca{i(\bt_{i+1})^\vee}{s_{i(\bt_i)}\cdots s_{i(\bt_1)}(-\pd)}<0.$$
Remark that this (and in fact the whole lemma) is valid for any
fundamental weight $\pd$ (the minuscule hypothesis implies more
precisely that this bracket has to be $-1$).
Let us calculate
$$\sca{\at_{i+1}^\vee}{-i(\pd)}= \sca{s_{\bt_1}\cdots
  s_{\bt_i}(\bt_{i+1})^\vee}{-i(\pd)}=
\sca{\bt_{i+1}^\vee}{s_{\bt_i}\cdots s_{\bt_1}(-i(\pd))}.$$
$$=\sca{i(\bt_{i+1})^\vee}{s_{i(\bt_i)}\cdots
  s_{i(\bt_1)}(-\pd)}<0.\ \ \ \ \ \ \ \ $$
that is to say $\sca{\at_{i+1}^\vee}{i(\pd)}>0$. Writing $\at_{i+1}$ in
terms of simple roots, we see that the coefficient of $\b$ has to be
strictly positive (in fact it has to be one because $\pd$ is
cominuscule). This exactly means that $\at_{i+1}\geq\b.$\hfill$\Box$

\vs 0.4 cm

It is now an easy check on the tables of \cite{bourb} to see that for
these roots and a minuscule weight $\pd$ we always have
$$\sca{\at_i^\vee}{\at_j}\geq0.$$
\vs -0.9 cm
\hspace{\stretch{10}}$\Box$

\begin{coro}
\label{debut}
With the above notations and the remark of lemma \ref{pgqb}, the fact
  that the expression $\phi=s_{i(\b_1)}\cdots s_{i(\b_n)}$ is reduced
  implies that  $\sca{i(\b_n)}{-\pd}<0$ or equivalently
  $\sca{\b_n}{i(\pd)}>0$. This is possible if and only if $\b_n=\b$.
\end{coro}

Let $k\in[1,n]$, if there exists an $i<k$ such that $\b_i=\b_k$
(resp. if there exists an $i>k$ such that $\b_i=\b_k$) we will denote
by $p(k)$ (resp. $n(k)$) the biggest (resp. smallest) integer
$i\in[1,k-1]$ (resp. $i\in[k+1,n]$) such that $\b_i=\b_k$.

\begin{coro}
\label{fin}
Let $j$ such that $\b_j=\b$.

{\rm (\i)} We have $\sca{\a_i^\vee}{\a_j}=0$ if for all $k\in[i+1,j]$,
$\sca{\b_i^\vee}{\b_k}=0$.

{\rm (\i\i)} Otherwise we have $\sca{\a_i^\vee}{\a_j}=1$ if $i>p(j)$
or if $i<p(j)$ and for all $k\in[i+1,p(j)]$,
$\sca{\b_i^\vee}{\b_k}=0$. In all other cases we have
$\sca{\a_i^\vee}{\a_j}=0$.
%
%the formula
%$$\sca{\a_i^\vee}{\a_n}=\left\{
%\begin{array}{cc}
%0&\ {\textit if}\ \ {\textit for}\ \ {\textit all}\ k\in[i-1,n],\
%\sca{\b_i^\vee}{\b_k}=0\\ 
%1&\ {\textit if}\ \b_k\neq\b\ {\textit for}\ k\in[i,n-1]\\ 
%0&\ {\textit otherwise}
%\end{array}\right..$$
\end{coro}

\dm
We have seen that
$\sca{\a_i^\vee}{\a_j}=\scal{\at_{n+1-i}^\vee}{\at_{n+1-j}}$ and
composing with $s_{\bt_1}\cdots s_{\bt_{n-j}}$ we can assume that
$n+1-j=1$ (ie $j=n$). We thus have to calculate
$\scal{\at_{n+1-i}^\vee}{\at_{1}}=\scal{\bt_1^\vee}{\at_{n+1-i}}=
\sca{\b^\vee}{\at_{n+1-i}}$ (we use here the fact that $R=R^\vee$ and
the fact that $\b_j=\b$). We first have to prove that
$\sca{\b^\vee}{\at_{n+1-i}}=0$ if for all $k\in[1,n+1-i]$,
$\scal{\bt_{n+1-i}^\vee}{\bt_k}=0$. 

And otherwise we have to prove that
$\sca{\b^\vee}{\at_{n+1-i}}=1$ if $n+1-i<n(1)$ or if $n+1-i>n(1)$ and
for all $k\in[n(1),n+1-i]$, $\scal{\bt_{n+1-i}^\vee}{\bt_k}=0$ and
that in all other cases we have $\sca{\b^\vee}{\at_{n+1-i}}=0$.

%$$\sca{\b^\vee}{\at_{n+1-i}}=\left\{
%\begin{array}{cc}
%1&\ {\textit if}\ \bt_k\neq\b\ {\textit for}\ k\in[2,n+1-i]\\ 
%0&\ {\textit otherwise}
%\end{array}\right..$$

(\i) In this case, it is easy to see that
$\at_{n+1-i}=\bt_{n+1-i}$ and we have the vanishing.

(\i\i) Let us define $\a=s_{\bt_2}\cdots s_{\bt_{n-i}}(\bt_{n+1-i})$. We have
$\at_{n+1-i}=s_{\bt_1}(\a)=s_\b(\a)$. And recall that the simple root
$\b$ always appears in $\at_{n+1-i}$ (lemma \ref{pgqb}) with
multiplicity 1 (because $\pd$ is a cominuscule weight). 

In the first case, we see that the simple root $\b$ does not appear
in $\a$. But we have $\at_{n+1-i}=s_\b(\a)=\a-\sca{\b^\vee}{\a}\b$
thus $\sca{\b^\vee}{\a}=-1$.

In the second case, applying lemma \ref{pgqb} to the sequence
$n(1),\cdots,n+1-i$ we see that the simple root $\b$ appears in
$s_{\bt_{n(1)}}\cdots s_{\bt_{n-i}}(\bt_{n+1-i})$ with multiplicity
1. As $\b$ does not appear in $\bt_2,\cdots,\bt_{n(1)-1}$, we see that
$\b$ appears in $\a$ with multiplicity 1. But we have
$\at_{n+1-i}=s_\b(\a)=\a-\sca{\b^\vee}{\a}\b$ thus
$\sca{\b^\vee}{\a}=0$.

We conclude because $\sca{\b^\vee}{\at_{n+1-i}}=
\sca{\b^\vee}{s_\b(\a)}=-\sca{\b^\vee}{\a}$.\hfill$\Box$

\begin{rema}
\label{fin2}
  The formula of corollary \ref{fin} is more simple if we use
  commutation relation beetween the simple root $\b_k$: let $j$ such
  that $\b_j=\b$, then we have $\sca{\a_i^\vee}{\a_j}=0$ if modulo
  commutation we can exchange $s_{\b_i}$ and $s_{\b_j}$. If not we
  also have $\sca{\a_i^\vee}{\a_j}=0$ if $i<p(j)$ and we can not
  commute $s_{\b_i}$ and $s_{\b_{p(j)}}$.
%works for $\sca{\a_i^\vee}{\a_j}$
%  if $\b_j=\b$ (the important point is that $\b$ appears in any root
%  with multiplicity at most one).
\end{rema}

Let us prove the following:

\begin{coro}
\label{rectif}
  We have the formula
$$\sum_{k=i+1,\ \b_k=\b}^n\sca{\a_i^\vee}{\a_k}=\left\{
  \begin{array}{cc}
1&\textit{if}\ \b_i\neq\b\\
0&\textit{if}\ \b_i=\b
  \end{array}\right..$$
\end{coro}

\dm
We apply the previous corollary. We know that $\b_n=\b$ and we can not
commute $s_{\b_i}$ and $s_{\b_n}$ (otherwise the expression would not
be reduced). Let $j$ be the smallest integer $k\in[i+1,n]$ such that
$\b_k=\b$ and we can not commute $s_{\b_i}$ and $s_{\b_k}$.

We have $\sca{\a_i^\vee}{\a_k}=0$ for all $k\in[i+1,n]$ with $\b_k=\b$
and $k\neq j$. For $k=j$, we have 
$$\sca{\a_i^\vee}{\a_k}=\left\{
  \begin{array}{cc}
1&\textit{if}\ \b_i\neq\b\\
0&\textit{if}\ \b_i=\b
  \end{array}\right..$$
\vs -0.9 cm\hfill$\Box$

\vs 0.4 cm

As is the proposition \ref{posi}, it is easy to check on the tables of
\cite{bourb} the following 

\begin{fact}
 If $\pd$ is minuscule and $(\a_i)_{i\in[1,n]}$ as above, then for all
  $i$ and $j$ in $[1,n]$, one has $\sca{\a_i^\vee}{\a_j}\leq2$
  with equality if and only if $\a_i=\a_j$.
\end{fact}

Let us prove the following corollary that we will need later:

\begin{coro}
\label{pgqmoins1}
  Let $i$, $x$ and $j$ in $[1,n]$. If $\sca{\a_i^\vee}{\a_x}=1$ then for
  all $j\in[1,n]$, one has
$$\sca{\a_i^\vee}{s_{\a_x}(\a_j)}\geq-1.$$
\end{coro}

\dm
We have
$$\sca{\a_i^\vee}{s_{\a_x}(\a_j)}= \sca{\a_i^\vee}{\a_j}-
\sca{\a_i^\vee}{\a_x}\sca{\a_x^\vee}{\a_j}=
\sca{\a_i^\vee}{\a_j}-\sca{\a_x^\vee}{\a_j}.$$
The preceding fact tells us that $\sca{\a_x^\vee}{\a_j}\leq2$ with
equality only if $\a_x=\a_j$. In case of equality we have
$\sca{\a_i^\vee}{\a_j}=\sca{\a_i^\vee}{\a_x}=1$ thus
$\sca{\a_i^\vee}{s_{\a_x}(\a_j)}=-1$. 

If $\a_x\neq\a_j$, then proposition \ref{posi} tells us that
$\sca{\a_i^\vee}{\a_j}\geq0$ and we have
$\sca{\a_x^\vee}{\a_j}\leq1$ thus
$\sca{\a_i^\vee}{s_{\a_x}(\a_j)}\geq-1$.\hfill$\Box$

\section{The Bott-Samelson resolutions}

In this section we briefly describe the Bott-Samelson construction which
give a resolution of any Schubert variety in $G/B$ and in $G/P$ for
any parabolic subgroup $P$. We describe this construction as
M. Demazure did in \cite{DE} we refer to this article for more
details.

\subsection{Construction}

Let $\phi\in W$ with $l(\phi)=n$. We recall in this paragraph M.
Demazure's construction \cite{DE} of a resolution of the dimension $n$
Schubert variety $X(\phi)=\overline{B\phi B/B}\subset G/B$
associated to a reduced decomposition $\phi=s_{\ga_{1}}\cdots
s_{\ga_{n}}$ with $\ga_i\in S$.

Let $w_0$ be the longuest element of $W$ and define the element
  $w=w_0\phi^{-1}w_0$. The preceding reduced expression leads to the
  reduced expression
$$w=s_{i(\ga_{n})}\cdots s_{i(\ga_{1})}.$$
If we choose any reduced expression
$$ww_0=s_{i(\ga_{n+1})}\cdots s_{i(\ga_{N})}$$
with $\ga_i\in S$ and $N=l(w_0)$, then $w_0=s_{i(\ga_{1})}\cdots
s_{i(\ga_N)}$ is a reduced expression of $w_0$. To keep the same
notation with \cite{DE}, let us note $\b_i=i(\ga_i)$, we have:
$$w_0=s_{\b_1}\cdots s_{\b_N},\ \ w=s_{\b_{n}}\cdots s_{\b_{1}}\ \ {\rm
  and}\ \ ww_0=s_{\b_{n+1}}\cdots s_{\b_{N}}.$$
With the sequence $(\b_i)_{i\in[1,N]}$, we define the following
sequence $(\a_i)_{i\in[1,N]}$ of roots by:
$$\a_1=\b_1,\ \a_2=s_{\b_1}(\b_2),\ \dots,\
\a_N=s_{\b_1}\cdots s_{\b_{N-1}}(\b_N).$$
The $\a_i$ are distincts and $\Del^+=\{\a_i\ /\ i\in[1,N]\}$. Define
$w_i=s_{\a_i}\in W$ (we will also for simplicity of notations sometimes
consider $w_i$ as an element of $G$). We have
$$w_i=s_{\b_1}\cdots s_{\b_{i-1}}s_{\b_i}s_{\b_{i-1}}\cdots
s_{\b_{1}},\ \ w=w_1\cdots w_{n},\ \ w_0=w_1\cdots w_N\ \ {\rm
  and}\ \ w_0^{-1}\phi=w_{N-n+1}\cdots w_N.$$ 
We define a sequence $(B_i)_{i\in[0,N]}$ of Borel subgroups containing
$T$ by induction: 
$$B_0=\Bt\ {\rm and}\ B_{i+1}=w_{i+1}(B_i).$$
Denote by $P_i$ the parabolic subgroup generated by $B_{i-1}$ and
$B_i$ we get a sequence of codimension one inclusions:
$$B_0\subset P_1\supset B_1\subset\cdots\supset B_{n-1}\subset
P_N\supset B_N.$$
Finaly we construct a sequence of varieties $(X_i)_{i\in[0,N]}$
endowed with a right action of $B_i$ by induction:
$$X_0=B_0\ {\rm and}\ X_{i+1}=X_i\times^{B_i}P_i$$
where the second term is the contracted product of $X_i$ and $P_i$
over $B_i$ (see \cite{DE} Par. 2.3.). The quotient $X_i/B_i$ is well
defined and we get a sequence of $\pu$-bundles $f_i$ with canonical
sections~$\sigma_i$:
$$X_0/B_0\xleftarrow{f_1}X_1/B_1\leftarrow\cdots\leftarrow
X_{N-1}/B_{N-1}\xleftarrow{f_N}X_N/B_N.$$
The scheme $X_i/B_i$ is 
the quotient of $P_1\times\cdots\times P_i$
by the right action of $B_1\times\cdots\times B_i$ given by
$$(p_1,\cdots,p_i)\cdot(b_1,\cdots,b_i)=(p_1b_1,\cdots,b_{i-1}^{-1}p_ib_i).$$
The projection $f_i$ sends the class of $(p_1,\cdots,p_i)$ to the
class of $(p_1,\cdots,p_{i-1})$ whereas the section $\sigma_i$ sends
the class of $(p_1,\cdots,p_{i-1})$ to the class of
$(p_1,\cdots,p_{i-1},w_i)$.

The multiplication morphism $P_1\times\cdots\times P_N\to G$ factorises
through $X_N\to G$ which is $P_N$ equivariant and in particular $B_N$
equivariant. We thus get a morphism
$$X_N/B_N\to G/B_N=G/B$$
which is birational and such that the restriction to
$\sigma_{N}\cdots\sigma_{n+1}(X_{n}/B_{n})$ is
birational on the Schubert variety $\overline{\Bt w_0^{-1}\phi
  B/B}\simeq X(\phi)$.
This construction gives us the resolution
$$\pi:X_n/B_n\to X(\phi).$$

Let $P$ be a parabolic subgroup containing $P$ and let $\bar\phi\in
W/W_P$. We want to construct a resolution of the Schubert variety
$X(\bar\phi)=\overline{B\bar\phi P/P}\subset G/P$.
For this choose $\phi$ the smallest element in the class
$\bar\phi$. The morphism $X(\phi)\to X(\bar\phi)$ induced by
the projection $G/B\to G/P$ is birational. So the morphism
$$\pi:X_n/B_n\to X(\bar\phi)$$
is a resolution. We will denote by $\Xt(\bar\phi)$ the scheme
$X_n/B_n$.

\begin{rema}
\label{echange}
  If we have $\sca{\b_i^\vee}{\b_{i+1}}=0$ for some $i$, then the
  Boot-Samelson resolution associated to the sequence
  $(\b_k)_{k\in[1,n]}$ is the same as the Boot-Samelson resolution
  associated to the sequence $(\b'_k)_{k\in[1,n]}$ where $\b'_k=\b_k$
  for $k\not\in\{i,i+1\}$, $\b'_i=\b_{i+1}$ and $\b'_{i+1}=\b_i$.
\end{rema}

\subsection{Curves and divisors on the Bott-Samelson resolution}

In his paper \cite{DE}, M. Demazure studies some special cycles on the
varieties $X_N/B_N$. Denote $Z_i=f_k^{-1}\cdots
f_{i+1}^{-1}(\im(\sigma_i))$. It is a divisor in $X_N/B_N$. For any
$K\subset[1,N]$ denote by 
$$Z_K=\bigcap_{i\in K}Z_i$$
which is a codimension $\vert K\vert$ subvariety of $X_N/B_N$. The
classes of the $Z_K$ form a basis of the Chow group of $X_N/B_N$
(cf. \cite{DE} Par. 4. prop. 1).
Remark that for any $k\in[1,N]$, we have $X_k/B_k=Z_{[k+1,N]}$. We can
in this way define subvarieties of $\Xt(\bar\phi)$: 
\begin{itemize}
\item denote by $D_i=Z_{\{i\}\cup[n+1,N]}$. This is a divisor on
  $\Xt(\bar\phi)$ and 
  these divisors form a basis of the Picard group of $\Xt(\bar\phi)$.
\item Define the curve $C_i=Z_{[1,N]-\{i\}}$. These curves for
  $i\in[1,n]$ form a basis of $A_1(\Xt(\bar\phi))$.
\end{itemize}

Denote by $\xi_K$ the class of $Z_K$ in the Chow group of
$X_N/B_N$. M. Demazure describes completely the Chow group of
$X_N/B_N$ in the following

\begin{theo}
{\rm (Demazure \cite{DE} Par. 4. prop. 1)}
  The Chow group of $X_N/B_N$ is generated over $\Z$ by
the $(\xi_i)_{i\in[1,N]}$ with the relations:
$$\xi_i\cdot\left(\sum_{j=1}^i\sca{\a_j^\vee}{\a_i}\xi_j\right)=0\ \ {\rm
    for}\ \ {\rm all}\ \ i\in[1,N].$$
\end{theo}

With the above notation we have
$\displaystyle{[C_i]=\prod_{j\neq i}\xi_j}$ and we can use the
previous theorem to prove

\begin{prop}
  We have 
$$[C_i]\cdot\xi_j=\left\{\begin{array}{cc}
0&\ {\rm for}\ i>j\\
1 &\ {\rm for}\ i=j\\
\sca{\b_i^\vee}{\b_j} &\ {\rm for}\ i<j
\end{array}\right..$$
\label{intersection}
\end{prop}

\dm
The preceding theorem leads by an easy induction to 
 
\begin{fact}
  We have the following formula in $A(X_N/B_N)$:
$$[C_i]\cdot\xi_j=\left\{\begin{array}{cc}
0&\ {\rm for}\ i>j\\
1 &\ {\rm for}\ i=j\\
\displaystyle{\sum_{k=1}^{j-i}(-1)^k\sum_{i=i_0<\cdots<i_k=j}
  \ \prod_{x=0}^{k-1}\sca{\a_x^\vee}{\a_{x+1}}}& \ {\rm for}\ i<j 
\end{array}\right..$$
\label{fact}
\end{fact}

We prove the following lemma to conclude the proof:

\begin{lemm} For $i<j$, we have the following formula:
$$\displaystyle{\sum_{k=1}^{j-i}(-1)^k\sum_{i=i_0<\cdots<i_k=j}
 \ \prod_{x=0}^{k-1}\sca{\a_x^\vee}{\a_{x+1}}}=\sca{\b_i^\vee}{\b_j}.$$
\label{lemmaform}
\end{lemm}

\dm
Let us first remark that the $\b_i$ can be constructed thanks to the
$\a_i$ in the following way:  
$$\b_1=\a_1,\ \b_2=s_{\a_1}(\a_2),\ \dots,\
\b_N=s_{\a_1}\cdots s_{\a_{N-1}}(\a_N).$$
Calculating
$$\sca{\b_i^\vee}{\b_j}=\sca{s_{\a_1}\cdots
  s_{\a_{i-1}}(\a_i)^\vee}{s_{\a_1}\cdots s_{\a_{j-1}}(\a_j)}
  =\sca{\a_i^\vee}{s_{\a_i}\cdots s_{\a_{j-1}}(\a_j)}$$ 
$$\ \!=-\sca{\a_i}{s_{\a_{i+1}}\cdots s_{\a_{j-1}}(\a_j)}.\ \ \ \ \ \ \ \
\ \ \ \ \ \ \ \ \ \ \ \ \ \ \ \ \ \ \ \ \ \ \ \ \ \ $$ 
Furthermore we can write
$$s_{\a_i}\cdots s_{\a_{j-1}}(\a_j)=\sum_{k=i}^jx_{k,j}\a_k$$
with $x_{k,j}\in\Z$ not depending on $i$. On the one hand, we get by
an easy induction the equality:
$$x_{i,j}=\displaystyle{\sum_{k=1}^{j-i}(-1)^k\sum_{i=i_0<\cdots<i_k=j}
  \prod_{x=0}^{k-1}\sca{\a_x^\vee}{\a_{x+1}}}.$$
On the other hand, we have
$$\sca{\a_i^\vee}{s_{\a_i}\cdots s_{\a_{j-1}}(\a_j)}
=\sum_{k=i}^jx_{k,j}\sca{\a_i^\vee}{\a_k}$$
and
$$-\sca{\a_i^\vee}{s_{\a_{i+1}}\cdots
s_{\a_{j-1}}(\a_j)}=-\sum_{k=i+1}^jx_{k,j}\sca{\a_i^\vee}{\a_k}$$
summing the two equalities we get
$$2\sca{\a_i^\vee}{s_{\a_i}\cdots
  s_{\a_{j-1}}(\a_j)}=\sca{\a_i^\vee}{s_{\a_i}\cdots
  s_{\a_{j-1}}(\a_j)}-\sca{\a_i^\vee}{s_{\a_{i+1}}\cdots
  s_{\a_{j-1}}(\a_j)}=x_{i,j}\sca{\a_i^\vee}{\a_i}$$
concluding the proof of the lemma.\hfill$\Box$

The proposition follows from fact \ref{fact} and lemma
\ref{lemmaform}. \hfill$\Box$

\begin{rema}
  The formulae of proposition \ref{intersection} are still valid on
  $\Xt(\bar\phi)$.
\end{rema}

Let us introduce some notations (see also \cite{DE}). If $\lambda$ is
a character of the Torus $T$ let us denote by $\L_i(\l)$ the
associated line bundle on $X_i/B_i$ (recall that $T\subset B_i$). Let
us now denote by $T_i$ the relative tangent sheaf of the
$\pu$-fibration $f_i:X_i/B_i\to X_{i-1}/B_{i-1}$. Thanks to \cite{DE}
Par. 2. Prop. 1. and an easy induction on $i$ we get the 

\begin{fact}
  Let us still denote $\L_i(\l)$ the correponding class in
  $A^*(X_i/B_i)$ then we have the formula:
$$\L_i(\l)=\sum_{k=1}^i\sca{\a_k^\vee}{\l}\cdot\xi_k.$$
\end{fact}

Furthermore, M. Demazure remarks (\cite{DE} Par. 2. remark following
Prop. 1.) that we have $T_i=\L_i(\a_i)$ so that we get the following

\begin{coro}
  \label{tgt}
  Let us still denote $T_i$ the correponding class in
  $A^*(X_i/B_i)$ then we have the formula:
  $$T_i=\sum_{k=1}^{i}\sca{\a_k^\vee}{\a_i}\cdot\xi_k.$$
\end{coro}

Remark that the factor of $\xi_i$ in $T_i$ is 2. We get the 

\begin{coro}
 Let $C$ be a curve on $X_i/B_i$. Suppose that for all $k\in[1,i]$ we
 have $[C]\cdot\xi_k\geq0$ and  $\sca{\a_k^\vee}{\a_i}\geq0$
 then for all $k$ we have
$$[C]\cdot (T_k-\xi_k)\geq0\ {\rm and}\ {\rm in}\
{\rm particular}\ [C]\cdot T_k\geq0$$
where we still denote by $T_k$ the pull-back of $T_k$ on $X_i/B_i$.
\end{coro}

Finally if $\phi$ is the smallest element in the class $\bar\phi\in
W/W_{P_\pd}$ with $\pd$ a minuscule weight, the results of the
proposition \ref{posi} gives us

\begin{coro}
\label{interposi}
  Let $C$ be a curve on $\Xt(\bar\phi)$ the resolution of
  $X(\bar\phi)$. Suppose that for all $k\in[1,n]$ we  have
  $[C]\cdot\xi_k\geq0$ then for all $k$ we have 
$$[C]\cdot (T_k-\xi_k)\geq0\ {\rm and}\ {\rm in}\
{\rm particular}\ [C]\cdot T_k\geq0$$
where we still denote by $T_k$ the pul-back of $T_k$ on $\Xt(\bar\phi)$.
\end{coro}

\begin{prop}
\label{intersectiontgt}
  We have
  $$[C_i]\cdot T_j=\left\{\begin{array}{cc}
0&\ {\rm for}\ i>j\\
\sca{\b_i^\vee}{\b_j} &\ {\rm for}\ i\leq j
\end{array}\right..$$
\end{prop}

\dm
Thanks to corollary \ref{tgt} the result is clear for $i>j$. Let
$i\leq j$ and let us use corollary \ref{tgt} and proposition
\ref{intersection} to get
$$[C_i]\cdot
T_j=\sum_{k=1}^j\sca{\a_k^\vee}{\a_j}[C_i]\cdot\xi_k=
\sum_{k=i+1}^{j-1}\sca{\a_k^\vee}{\a_j}
\cdot\sca{\b_i^\vee}{\b_k}+\sca{\a_i^\vee}{\a_j}
+2\sca{\b_i^\vee}{\b_j}.$$ 

\begin{lemm}
  We have the formula
$$\sum_{k=i+1}^{j-1}\sca{\a_k^\vee}{\a_j}
\cdot\sca{\b_i^\vee}{\b_k}=-\sca{\a_i^\vee}{\a_j}
-\sca{\b_i^\vee}{\b_j}$$
\label{formule2}
\end{lemm}

\dm
Because the construction of $\a_i$ in terms of $\b_i$ is symetric
to the construction of $\b_i$ in terms of $\a_i$ the formula of lemma
\ref{lemmaform} is valid when we exchange the roles of the $\a_i$ and
of the $\b_i$ so we get for $k<j$:
$$\sca{\a_k^\vee}{\a_j}=\displaystyle{\sum_{u=1}^{j-k}(-1)^u
  \sum_{k=i_0<\cdots<i_u=j}\
  \prod_{x=0}^{u-1}\sca{\b_{i_x}^\vee}{\b_{i_{x+1}}}}.$$ 
We thus obtain
$$\sum_{k=i+1}^{j-1}\sca{\a_k^\vee}{\a_j}
\cdot\sca{\b_i^\vee}{\b_k}=\sum_{k=i+1}^j
\displaystyle{\sum_{u=1}^{j-k}(-1)^u \sum_{k=i_0<\cdots<i_u=j}\ 
  \prod_{x=0}^{u-1}\sca{\b_{i_x}^\vee}{\b_{i_{x+1}}}}
\cdot\sca{\b_i^\vee}{\b_k}.$$ 
If we set $i_{-1}=i$ we get
$$\sum_{k=i+1}^{j-1}\sca{\a_k^\vee}{\a_j}
\cdot\sca{\b_i^\vee}{\b_k}=\sum_{u=1}^{j-i-1}
\displaystyle{\sum_{k=i+1}^{j-u}(-1)^u 
  \sum_{i=i_{-1}<k=i_0<\cdots<i_u=j}\
  \prod_{x=-1}^{u-1}\sca{\b_{i_x}^\vee}{\b_{i_{x+1}}}}\ \ \ \ \ \ \ \
\ \ \ \ \ \ \ \ $$
$$\ \ \ \ \ \ \ =\sum_{u=1}^{j-i-1}
\displaystyle{(-1)^u 
  \sum_{i=i_{-1}<i_0<\cdots<i_u=j}\
  \prod_{x=-1}^{u-1}\sca{\b_{i_x}^\vee}{\b_{i_{x+1}}}}.$$
$$\ \!=\sum_{u=2}^{j-i}
\displaystyle{(-1)^{u+1} 
  \sum_{i=i_0<\cdots<i_u=j}\
  \prod_{x=0}^{u}\sca{\b_{i_x}^\vee}{\b_{i_{x+1}}}}.$$
$$\ \ \ \ \ \ \ \ \ \ \ \ \ \ \ \ \ \ =-
\sca{\b_i^\vee}{\b_j}+\sum_{u=1}^{j-i}
\displaystyle{(-1)^{u+1} 
  \sum_{i=i_0<\cdots<i_u=j}\
  \prod_{x=0}^{u}\sca{\b_{i_x}^\vee}{\b_{i_{x+1}}}}.$$
$$\ \ =
-\sca{\b_i^\vee}{\b_j}-\sca{\a_i^\vee}{\a_j}.\ \ \ \ \ \ \ \ \ \ \ \ \ \
\ \ \ \ \ \ \ \ \ \ \ \ $$ 
\vs -0.9 cm\hfill$\Box$

\vs 0.4 cm

This lemma with the preceding formula ends the proof.\hfill$\Box$

\section{Some more curves on $\Xt(\bar\phi)$}

\subsection{Effective and contracted curves}

In this paragraph, we study some more curves on $\Xt(\bar\phi)$. In
particular those who are contracted by the projection
$\pi:\Xt(\bar\phi)\to X(\bar\phi)$.

Let us look at the restriction of $\pi$ on the curve $C_j$. M. Demazure
(\cite{DE} Par. 3 Theorem 1) proves that the curve is contracted if
and only if $l(w_1\cdots w_{j-1}w_{j+1}\cdots w_n)>n$. But a simple
calculation gives
$$w_1\cdots w_{j-1}w_{j+1}\cdots w_n=s_{\b_j}w$$
so the curve is not contracted in $G/B$ if and only if
$l(s_{\b_j}w)=l(w)-1$ in $W$ and is not contracted in $G/P$ if this
equality is true for the minimal representatives in $W$ of
$\overline{s_{\b_j}w}$ and $\bar w$ in $W/W_P$. This means that there
exists a minimal reduced expression of $\bar w$ beginning with
$s_{\b_j}$. But for any reduced expression $\bar w=s_{\b_n}\cdots
s_{\b_1}$ we have seen in corollary \ref{debut} that we must have
$\b_n=\b$ (where $\b$ is the unique simple root such that
$\sca{\b^\vee}{i(\pd)}=1$). This would imply that $\b_j=\b=\b_n$. In
the other cases the curve $C_j$ is contracted (in general i.e. when
$\pd$ is not minuscule and even not fundamental, the curve $C_j$ is
contracted if and only if $\scal{\b_j^\vee}{i(\pd)}>0$). 

Let us consider the case $\b_j=\b=\b_n$. Define $n(j)$ the smallest
integer $k$ such that $k>j$ and $\b_k=\b_j$. In this case the curves
$C_j$, $C_{n(j)}$ and $C_n$ are not contracted. The morphism $\pi$
induces an isomorphism (because these curves are $\pu$) onto their
images which are respectively the Schubert varieties associated to
$s_{\b_j}$, $s_{\b_{n(j)}}$ and $s_{\b_n}$ (see \cite{DE}).
We see that they have the same image.

Let $t$ be a point in the commun image of the curve $C_j$ and
$C_{n(j)}$. The antecedent of the point $t$ in $C_j$ (resp. $C_{n(j)}$) is
the image in $\Xt(\bar\phi)$ of a $n$-uple
$(w_1,\cdots,w_{j-1},x(t),w_{j+1},\cdots,w_n)\in P_1\times\cdots\times
P_n$ (resp. $(w_1,\cdots,w_{n(j)-1},y(t),w_{n(j)+1},\cdots,w_n)$). We
thus have the equation
$$w_1\cdots w_{j-1}\cdot x(t)\cdot w_{j+1}\cdots w_n=w_1\cdots
w_{n(j)-1}\cdot y(t)\cdot w_{n(j)+1}\cdots w_n.$$
If we consider the curve $\Ct_j$ parametrized by $t$ defined by the
images of
$$(w_1,\cdots,w_{j-1},w_{j}^{-1}x(t),w_{j+1},
\cdots,w_{n(j)-1},w_{n(j)}y(t)^{-1},w_{n(j)+1},\cdots,w_n)$$ 
in $\Xt(\bar\phi)$ we see that its image by $\pi$ is $w_1\cdots
w_{j-1}w_{j+1} \cdots w_{n(j)-1}w_{n(j)+1}\cdots w_{n-1}$ a
constant. The curve $\Ct_j$ is contracted by $\pi$.

\begin{lemm}
  We have $[\Ct_j]=[C_j]-[C_{n(j)}]$.
\end{lemm}

\dm
The projection of $\Ct_j$ and $C_j$ on $X_{n(j)-1}/B_{n(j)-1}$ are the
same. This implies that $[C_j]-[\Ct_j]=a[C_{n(j)}]$ with
$a\in\Z$. Apply $\pi_*$ to this equation to get
$\pi_*[C_j]-\pi_*[\Ct_j]=a\pi_*[C_{n(j)}]$. But we have
$\pi_*[C_j]=\pi_*[C_{n(j)}]$ and $\pi_*[\Ct_j]=0$ thus
$a=1$.\hfill$\Box$

\begin{prop}
  The classes $[\Ct_j]$ generate $A_1(\Xt(\bar\phi))$ over
  $\Z$. Furthermore they generate the cone of effective curves
  i.e. they generate the extremal rays.
\end{prop}

\dm
The first assertion is trivial because the classes $[C_j]$ generate
$A_1(\Xt(\bar\phi))$ over $\Z$. 

For the second, we proceed by induction on $j$: we prove that the
classes $[\Ct_k]$ for $k\leq j$ generated the effective cone of
$X_j/B_j$ (by abuse of notation we still denote by $[\Ct_k]$ the image
of the class $[\Ct_k]$ in $X_j/B_j$). It is true for $j=1$ assume it
is true for $j-1$ and let
$$[C]=\sum_{k=1}^ja_k[\Ct_k]$$
the class of an effective curve. By projection on $X_{j-1}/B_{j-1}$ we
obtain the class
$${f_j}_*[C]=\sum_{k=1}^{j-1}a_k[\Ct_k]$$
which has to be effective so by induction we have $a_k\geq0$ for
$k<j$. Now by projection on $G/B$ we get 
$$\pi_*[C]=\sum_{k=1}^{j-1}a_k\pi_*[\Ct_k].$$
The class $[\Ct_j]$ is not contracted. The only classes $[\Ct_k]$ that
are not contracted by $\pi$ are such that $[\Ct_k]=[C_k]$ and
$l(s_{\b_k}w)=l(w)-1$. The image is then the Schubert variety
associated to $s_{\b_k}$. The first condition implies that for these
not contracted curves, all the $\b_k$ are distinct. But the associated
Schubert varieties are independent in $A_1(G/B)$ and because the image
is effective we have $a_k\geq0$ for all those $k$ and in particular
$a_j\geq0$.\hfill$\Box$

\subsection{Curves on contracted divisors}

Let $x\in[1,n]$ such that the divisor $D_x$ is contracted by $\pi$. We
are going to construct special curves on $D_x$ (recall that $[D_x]=\xi_x$). 

\begin{lemm}
  There exists $i\in[1,n]$ such that $[C_i]\cdot\xi_x=-1$.
\label{existmoins1}
\end{lemm}

\dm
Recall that we have (proposition \ref{intersection})
$$[C_i]\cdot\xi_x=\left\{\begin{array}{cc}
0&\ {\rm for}\ i>x\\
1 &\ {\rm for}\ i=x\\
\scal{\b_i^\vee}{\b_x} &\ {\rm for}\ i<x
\end{array}\right..$$
We have to choose $i<x$ and for such an $i$, as the group is simply
laced we have $[C_i]\cdot\xi_x=-1$, $0$ or $2$. If for all $i<x$ this
intersection is zero then for all $i<x$ the symetry $s_{\b_i}$
commutes with $s_{\b_x}$ so that the reduced expression
$w=s_{\b_n}\cdots s_{\b_1}$ can be written $w=s_{\b_n}\cdots
s_{\b_{x+1}}s_{\b_{x-1}}\cdots s_{\b_1}s_{\b_x}$. We have a reduced
expression 
$$\bar\phi=s_{\ga_x}s_{\ga_1}\cdots s_{\ga_{x-1}}s_{\ga_{x+1}}\cdots
s_{\ga_n}$$ 
meaning that the image of $D_x$ in $X(\bar\phi)$ is a moving
divisor. This is impossible because $D_x$ is contracted. Let $i$ be
the biggest $i<x$ such that $[C_i]\cdot\xi_x\neq0$. If the
intersection is 2 this means that $\b_i=\b_x$. But because for all
$k\in[i+1,x-1]$, we have $\scal{\b_k^\vee}{\b_x}=0$, we see that $s_{\b_x}$
commutes with all $s_{\b_k}$ with $k\in[i+1,x-1]$. We have:
$$\bar\phi=s_{\ga_1}\cdots s_{\ga_{i-1}}s_{\ga_{i}}
s_{\ga_{i+1}}\cdots s_{\ga_{x-1}}s_{\ga_{x}}s_{\ga_{x+1}} \cdots
s_{\ga_n}= s_{\ga_1}\cdots
s_{\ga_{i-1}}s_{\ga_{i}}s_{\ga_{x}}s_{\ga_{i+1}}\cdots
s_{\ga_{x-1}}s_{\ga_{x+1}}\cdots s_{\ga_n}$$
$$=s_{\ga_1}\cdots s_{\ga_{i-1}}s_{\ga_{i+1}}\cdots
s_{\ga_{x-1}}s_{\ga_{x+1}}\cdots s_{\ga_n}\ \ \ \ \ \ \ \ \ \ \ \ \ \
\ \ \ \ \ \ \ \ \ \ \ \ \ \ \ \ \ \ \ \ \ \ \ \ \ \ \ \ \ \ \ \ \ \ \
\ \ \ \ \ \ \ \ \ \ \ \ \ \ $$ 
that is to say the expression $\bar\phi=s_{\ga_1}\cdots
s_{\ga_{n}}$ was not reduced, a contradiction.\hfill$\Box$

\begin{rema}
\label{existence}
  In particular there exists an  $i\in[1,n]$ such that
  $\sca{\a_i^\vee}{\a_x}=1$ (choose the $i$ of the preceding proof and
  we have $\sca{\a_i^\vee}{\a_x}=-\sca{\b_i^\vee}{\b_x}=1$).
\end{rema}

Let $i\in[1,n]$ and let us define the following classes of curves:
$$[\Ch_i]=[C_i]+\sum_{k=i+1}^n\sca{\a_i^\vee}{\a_k}[C_k].$$

\begin{lemm}
  We have the formulae
$$[\Ch_i]\cdot \xi_j=\delta_{i,j}\ \ {\textit and}\ \ [\Ch_i]\cdot
T_j=\left\{\begin{array}{cc}
0&\ {\rm for}\ i>j\\
\sca{\a_i^\vee}{\a_j} &\ {\rm for}\ i\leq j
\end{array}\right..$$
\label{interchap}
\end{lemm}

\dm
We use proposition \ref{intersection} and lemma \ref{formule2} to get
$$[\Ch_i]\cdot \xi_j=\left([C_i]+\sum_{k=i+1}^n
  \sca{\a_i^\vee}{\a_k}[C_k]\right)\cdot\xi_j\ \ \ \ \ \ \ \ \ \ \ \ \
\ \ \ \ \ \ \ \ \ \ \ \ $$
$$\ \ \ \ \ \ \ \ \ \ \ \ \ \ \ \ \ \ \ \ \ \ \!\ =\left\{\begin{array}{cc}
0&\ {\rm for}\ i>j\\
1&\ {\rm for}\ i=j\\
\displaystyle{\sca{\b_i^\vee}{\b_j}+\sum_{k=i+1}^{j-1}
  \sca{\a_i^\vee}{\a_k}\sca{\b_k^\vee}{\b_j}+\sca{\a_i^\vee}{\a_j}}
&\ {\rm for}\ i<j
\end{array}\right.$$
$$\ \ \ \ \ \ \ \ \ \ \ \ \ \ \ \ \ \ =\left\{\begin{array}{cc}
0&\ {\rm for}\ i>j\\
1&\ {\rm for}\ i=j\\
\displaystyle{\sca{\b_i^\vee}{\b_j}-\sca{\b_i^\vee}{\b_j}-
  \sca{\a_i^\vee}{\a_j}+\sca{\a_i^\vee}{\a_j}}&\  
  {\rm for}\ i<j
\end{array}\right.$$
proving the first formula. For the second one we use proposition
\ref{intersectiontgt} and lemma \ref{formule2} to get 
$$[\Ch_i]\cdot T_j=\left([C_i]+\sum_{k=i+1}^n
  \sca{\a_i^\vee}{\a_k}[C_k]\right)\cdot T_j\ \ \ \ \ \ \ \ \ \ \ \ \
\ \ \ \ \ \ \ \ \ \ \ \ $$
$$\ \ \ \ \ \ \ \ \ =\left\{\begin{array}{cc}
0&\ {\rm for}\ i>j\\
\displaystyle{\sca{\b_i^\vee}{\b_j}+\sum_{k=i+1}^j
  \sca{\a_i^\vee}{\a_k}\sca{\b_k^\vee}{\b_j}} &\ {\rm for}\ i\leq j
\end{array}\right.$$
$$\ \ \ \ \ \ \ \ \ \ \ \ \ \ \ \ \ \ \ \ \ =\left\{\begin{array}{cc}
0&\ {\rm for}\ i>j\\
1&\ {\rm for}\ i=j\\
\displaystyle{\sca{\b_i^\vee}{\b_j}-\sca{\b_i^\vee}{\b_j}-
  \sca{\a_i^\vee}{\a_j}+2\sca{\a_i^\vee}{\a_j}}&\  
  {\rm for}\ i<j
\end{array}\right.$$
concuding the proof.\hfill$\Box$

\vs 0.4 cm

Now let $i\in[1,n]$ such that $\sca{\a_i^\vee}{\a_x}=1$ (there exists
such an $i$ thanks to remark \ref{existence}). We define the class:
$$[\Ga_{x,i}]=[\Ch_i]-\sca{\a_i^\vee}{\a_x}[\Ch_x]=[\Ch_i]-[\Ch_x]$$
and prove the following:

\begin{prop}
\label{recouvre}
  We have:

{\rm (\i)} $[\Ga_{x,i}]\cdot\xi_x=-1$ so all curves $C\in[\Ga_{x,i}]$ are
contained in $D_x$.

{\rm (\i\i)} The scheme $\Mor{[\Ga_{x,i}]}{\Xt(\bar\phi)}$ is
irreducible and smooth (in particular non empty).
 
{\rm (\i\i\i)} The open part $\displaystyle{D_x-\bigcup_{k\neq
    x}(D_x\cap D_k)}$ of the divisor $D_x$ is covered by curves
    $C\in[\Ga_{x,i}]$.

{\rm (\i v)} All curves $C\in[\Ga_{x,i}]$ are contracted by $\pi$.
\end{prop}

\dm
(\i) This is a simple application of lemma \ref{interchap}.

\vs 0.2 cm

(\i\i) Recall that $\Xt(\bar\phi)$ is a sequence of $\pu$-bundles. We
proceed by induction on the $X_j/B_j$ (by abuse of notation, we still
denote by $[\Ga_{x,i}]$ the push-forward of $[\Ga_{x,i}]$ in $A_1(X_j/B_j)$). 
Let us denote by $\vp:X\to Y$ the morphism $f_j:X_j/B_j\to
X_{j-1}/B_{j-1}$ and by $T$ the relative tangent sheaf. We have a
section $\sigma=\sigma_j$ of $\vp$ and we denote by $\xi=\xi_j$ the
divisor image of the section. We have:
$$\sigma_*\vp_*[\Ga_{x,i}]= \left\{
\begin{array}{cc}
0& \ {\rm
 for}\ j\leq i\\{}
[\Ga_{x,i}]-\sca{\a_i^\vee}{\a_j}\cdot[C_j] & \ {\rm
 for}\ i<j<x\\{}
[\Ga_{x,i}]&\ {\rm for}\ j=x \\{} 
[\Ga_{x,i}]-\sca{\a_i^\vee}{s_{\a_x}(\a_j)}\cdot[C_j] & \ {\rm
 for}\ j>x
\end{array}\right..$$
Proposition \ref{intersection} and lemma \ref{interchap} give us
$$[\Ga_{x,i}]\cdot
\xi=\left\{\begin{array}{cc}
1&\ {\rm for}\ j=i\\
-1&\ {\rm for}\ j=x\\
0&\ {\rm otherwise}
\end{array}\right.,\ \ \ \ \ \sigma_*\vp_*[\Ga_{x,i}]\cdot
\xi=\left\{\begin{array}{cc}
0& \ {\rm
 for}\ j\leq i\\{}
-\sca{\a_i^\vee}{\a_j} & \ {\rm
 for}\ i<j\leq x\\{}
-\sca{\a_i^\vee}{s_{\a_x}(\a_j)} & \ {\rm
 for}\ j>x
\end{array}\right.$$
and
$$[\Ga_{x,i}]\cdot T=\left\{\begin{array}{cc}
0&\ {\rm for}\ j<i\\
\sca{\a_i^\vee}{\a_j} &\ {\rm for}\ i\leq j<x\\
\sca{\a_i^\vee}{s_{\a_x}(\a_j)}&\ {\rm for}\ j\geq x
\end{array}\right..$$
Let us denote by $[\Ga]$ the class of $[\Ga_{x,i}]$ in $X=X_j/B_j$ and
let $f\in\Mor{\vp_*[\Ga]}{Y}$. We want to study the fiber over $f$ of
the morphism
$$\Mor{[\Ga]}{X}\to\Mor{\vp_*[\Ga]}{Y}$$ 
that is to say the
morphisms $f'\in\Mor{[\Ga]}{X}$ such that $f=\vp\circ f'$. We look
for a section of the $\pu$-bundle $\vp$ pulled-back by $f$. Let $E$ be the
rank two vector bundle defining the $\pu$-bundle. We can choose $E$
such that $f^*E=\opu\oplus\opu(a)$ with $a\geq0$. 

The section $f\circ\sigma$
is given by a surjection $f^*E\to\opu(z)$ with
$2z-a=\sigma_*\vp_*[\Ga]\cdot \xi$.
A morphism $f'$ is simply given by a surjection $f^*E\to\opu(y)$ such
that $y+z-a=[\Ga]\cdot \xi$ and $2y-a=[\Ga]\cdot T$.

\begin{rema}
\label{xyz}
  The section $f\circ\sigma$ always exists. We must thus have $z=0$ or
  $z\geq a$. This implies that
  \begin{itemize}
  \item if $\sigma_*\vp_*[\Ga]\cdot \xi=0$ then $a=z=0$ ; 
  \item if $\sigma_*\vp_*[\Ga]\cdot \xi<0$ then $z=0$ and
    $a=-\sigma_*\vp_*[\Ga]\cdot \xi$ ;
  \item if $\sigma_*\vp_*[\Ga]\cdot
    \xi>0$ then $\sca{\a_i^\vee}{s_{\a_x}(\a_j)}=-\sigma_*\vp_*[\Ga]\cdot
    \xi<0$ and in fact $\sca{\a_i^\vee}{s_{\a_x}(\a_j)}=-1$ (corollary
    \ref{pgqmoins1}). In this case we have $2z-a=1$ and this implies
    $z=a=1$.
  \end{itemize}

The section $f'$ will exist if there exists an integer $y$ such that
  $y=0$ or $y\geq a$. In the case $j=i$, we have $z=a=0$, $y=1$ and $f'$
  exists. In the case $j=x$ we have $y=z=0$, $a=1$ and $f'$ exists. In
  the other cases we always have $y+z-a=0$. 
This implies that
  \begin{itemize}
  \item if $\sigma_*\vp_*[\Ga]\cdot \xi=0$ then $y=0$ ; 
  \item if $\sigma_*\vp_*[\Ga]\cdot \xi<0$ then $y=a$ ;
  \item if $\sigma_*\vp_*[\Ga]\cdot \xi>0$ then $y=0$.
  \end{itemize}
In conclusion there always exists a section $f'$ of $f$ with the
required invariants.
\end{rema}

We
will use the following proposition (see \cite{PE} Prop. 4):

\begin{prop}
  Let $\vp:X\to Y$ a $\pu$-bundle with relative tangent sheaf $T$ and
  let $[\Ga]\in A_1(X)$ such that $[\Ga]\cdot T\geq0$, then
  $\Mor{[\Ga]}{X}$ is an open subset of a projective bundle over
  $\Mor{\vp_*[\Ga]}{Y}$. In particular, if $\Mor{\vp_*[\Ga]}{Y}$ is
  irreducible, the same is true for $\Mor{[\Ga]}{X}$ as soon as it is
  non empty.
\label{pe}
\end{prop}

This proposition with be usefull for the fibration $f_j$ if we have
$[\Ga_{x,i}]\cdot T_j\geq0$. The only cases where the previous proposition
does not apply is when $\sca{\a_i^\vee}{s_{\a_x}(\a_j)}<0$ and in fact
$\sca{\a_i^\vee}{s_{\a_x}(\a_j)}=-1$ (lemma \ref{pgqmoins1}). There
are two distinct cases where this may occur. If $j=x$ then
$[\Ga_{x,i}]\cdot \xi_j=-1$ and $[\Ga_{x,i}]\cdot T_j=-1$. If $j>x$ and
$\sca{\a_i^\vee}{s_{\a_x}(\a_j)}=-1$ then $[\Ga_{x,i}]\cdot \xi_j=0$ and
$[\Ga_{x,i}]\cdot T_j=-1$.

The first case $j=x$ is treated thanks to the

\begin{lemm}
\label{sectionpaspb}
  Let $\vp:X\to Y$ a $\pu$-bundle  with relative tangent sheaf $T$ and
  with a section $\sigma$. Denote $\xi$ the divisor $\sigma(Y)$ and let
  $[\Ga]\in A_1(X)$ such that $[\Ga]\cdot \xi=-1$, $[\Ga]\cdot T=-1$
  and $\sigma_*\vp_*[\Ga]\cdot \xi=-1$. Suppose that
  $\Mor{\vp_*[\Ga]}{Y}$ is normal then we have
$$\Mor{[\Ga]}{X}\simeq\Mor{\vp_*[\Ga]}{Y}.$$
\end{lemm}

\dm
Let $f\in\Mor{\vp_*[\Ga]}{Y}$, we have to prove (by Zariski Main
theorem) that there is exactly
one morphism $f'\in\Mor{[\Ga]}{X}$ such that $f=\vp\circ f'$. But with
the above notation and thanks to remark \ref{xyz} we have $y=z=0$ and
$a=1$. The morphism $f'$ has to be $\sigma\circ f$.\hfill$\Box$

\vs 0.4 cm

The second case $j>x$ is treated thanks to the

\begin{lemm}
\label{sectionpb}
  Let $\vp:X\to Y$ a $\pu$-bundle  with relative tangent sheaf $T$ and
  with a section $\sigma$. Denote $\xi$ the divisor $\sigma(Y)$ and let
  $[\Ga]\in A_1(X)$ such that $[\Ga]\cdot \xi=0$, $[\Ga]\cdot T=-1$
  and $\sigma_*\vp_*[\Ga]\cdot \xi=1$. Suppose that
  $\Mor{\vp_*[\Ga]}{Y}$ is normal then we have
$$\Mor{[\Ga]}{X}\simeq\Mor{\vp_*[\Ga]}{Y}.$$
\end{lemm}

\dm
Let $f\in\Mor{\vp_*[\Ga]}{Y}$, we have to prove (by Zariski Main
theorem) that there is exactly
one morphism $f'\in\Mor{[\Ga]}{X}$ such that $f=\vp\circ f'$. But with
the above notation and thanks to remark \ref{xyz} we have $y=0$ and
$z=a=1$.
The morphism $f'$ is given by the unique self-negative section of
$\p_{\pu}(f^*E)$.\hfill$\Box$

\vs 0.2 cm

(\i\i\i) Let us note that thanks to remark \ref{xyz}, lemma
\ref{sectionpaspb} and lemma \ref{sectionpb} there are curves
$C\in[\Ga_{x,i}]$ such that $C$ is not contained in  any intersection
$D_x\cap D_j$ (we always have $C\subset D_x$) and thus always meet
the open part $\displaystyle{D_x-\bigcup_{k\neq x}(D_x\cap D_k)}$
of the divisor $D_x$. 

But the orbit of the unipotent part $U$ of $B$ acting on $D_x$ is
exactly $\displaystyle{D_x-\bigcup_{k\neq x}(D_x\cap D_k)}$. Translating
$C$ thanks to the action of $U$ we see that the curves $C\in[\Ga_{x,i}]$
cover $\displaystyle{D_x-\bigcup_{k\neq x}(D_x\cap D_k)}$. 

(\i v) We have seen that all the curves $[\Ct_k]$ are contracted by
$\pi$ except $[\Ct_n]$. We just have to prove that the coefficient $a_n$
of $[\Ct_n]$ in $[\Ga_{x,i}]$ is zero. Let us set
$$A=\sum_{k=i+1,\ \b_k=\b}^n\sca{\a_i^\vee}{\a_k}-\sum_{k=x+1,\
  \b_k=\b}^n\sca{\a_x^\vee}{\a_k}.$$
We have
$$a_n=\left\{
\begin{array}{cc}
A& \ {\rm
 if}\ \b_i\neq\b\ {\rm and}\ \b_x\neq\b \\{}
A+1 &  \ {\rm
 if}\ \b_i=\b\ {\rm and}\ \b_x\neq\b \\{}
A-1 &  \ {\rm
 if}\ \b_i\neq\b\ {\rm and}\ \b_x=\b \\{}
A & \ {\rm
 if}\ \b_i=\b\ {\rm and}\ \b_x=\b
\end{array}\right..$$
We now apply corollary \ref{rectif} to see that $a_n=0$ in all
cases.\hfill$\Box$ 

\begin{rema}
  If the fiber of the projection $\pi:D_x\to \pi(D_x)$ is a curve then
  its class has to be $[\Ga_{x,i}]$. In general, the generic fiber is
  covered by curves in the class $[\Ga_{x,i}]$. For more details on the
  fiber of the Bott-Samelson resolution see {\rm \cite{Gaussent}}.
\end{rema}

\section{The scheme of morhisms for $\Xt(\bar\phi)$}

\subsection{Irreducibility}

We will prove in this paragraph that for some classes $\a\in
A_1(\Xt(\bar\phi))$ the scheme of morphisms $\Mor{\a}{\Xt(\bar\phi)}$
is irreducible and smooth. We will essentially need proposition
\ref{pe} (see \cite{PE} Prop. 4).

Let us now consider a class $\a\in A_1(\Xt(\bar\phi))$ such that
$\a\cdot \xi_i\geq0$ for all $i\in[1,n]$. Thanks to corollary
\ref{interposi} we know that $\a\cdot T_i\geq0$ and $\a\cdot
(T_i-\xi_i)\geq0$.

\begin{prop}
\label{irredenhaut}
{\rm (\i)}  The scheme of morphisms $\Mor{\a}{\Xt(\bar\phi)}$ is
irreducible and smooth of dimension
$$\int_\a c_1\big(T_{\Xt(\bar\phi)}\big)+\dim\big(\Xt(\bar\phi)\big).$$

{\rm (\i\i)} If the class $\a$ is such that $\a\cdot\xi_x=0$ for all
$x\in[1,n]$ with $D_x$ a contracted divisor, then a general element
$f\in\Mor{\a}{\Xt(\bar\phi)}$ is contained in the regular locus of
$\pi$. 
\end{prop}

\vs 0.2 cm

\dm
(\i) We proceed by induction, for the first step, we have to study the
scheme of morphisms from $\pu$ to $\pu$. This scheme is
irreducible and smooth. We go by induction thanks to proposition
\ref{pe}. We only have to prove that the scheme is non empty. However
with the notations of the preceding paragraph for $\pu$-fibrations, we
have $f^*E=\opu\oplus\opu(a)$ with section $\sigma$ given by a
surjection $f^*E\to\opu(z)$ and we look for a section
$f^*E\to\opu(y)$. Because of the relations $\a\cdot T_i\geq0$ and $\a\cdot
(T_i-\xi_i)\geq0$ we see that $y\geq z$ and $y\geq a-z$. This implies
that $y\geq a$ proving the existence of a surjection
$f^*E\to\opu(y)$.

(\i\i) Let $f$ a general element. Thanks to the discussion above, we
may assume that this element will meet the non contracted divisors
$D_i$ in distinct points and will not meet the contracted
divisors. In particular $f$ will never meet intersections
$D_i\cap D_j$ with $i\neq j$. In particular, the only $B$-orbits of
the Bott-Samelson resolution that $f$ will meet are the dense orbit
and the orbits dense in $D_i$ for a non contracted divisor. These
orbits are contained in the regular locus so this in particular proves
that $f$ is contained in the regular locus.\hfill$\Box$

\subsection{Smoothing curves on $\Xt(\bar\phi)$}

Let $\a\in A_1(\Xt(\bar\phi))$ as above.

\begin{lemm}
  There exists $\ft\in\Mor{\a}{\Xt(\bar\phi)}$ such that $\ft(\pu)$
  in not contained in any $D_i$ and does
  not meet any intersection $D_i\cap D_j$.
\end{lemm}

\dm
Because the scheme $\Mor{\a}{\Xt(\bar\phi)}$ is irreducible, if it
exists, a general morphism will have the required property.

Let $i<j$, we construct this curve $\ft$ by induction on the
$\pu$-fibrations. For all fibrations except for the fibrations $f_i$
and $f_j$, we take any section. 

For the fibration $f_i$ we have by induction a morphism
$\ft_{i-1}:\pu\to X_{i-1}/B_{i-1}$ and with the notations of the proof of
the previous proposition: a rank 2 vector bundle
${\ft_{i-1}}^*E=\opu\oplus\opu(a)$ with $a\geq0$ ; a surjection
$\ft_{i-1}^*E\to\opu(z)$ (corresponding to the divisor $D_i$) and we
look for a surjection $\ft_{i-1}^*E\to\opu(y)$. With our hypothesis on
$\a$ we have $y\geq z$ and $y\geq a-z$ (cf. proof of the preceding
proposition) so there always exists a section and we can choose it
such that the image is not contained in $D_i$ (because $y\geq
a-z$).

For the fibration $f_j$ we have by induction a morphism
$\ft_{j-1}:\pu\to X_{j-1}/B_{j-1}$ and with the notations of the proof of
the previous proposition: a rank 2 vector bundle
$\ft_{j-1}^*E=\opu\oplus\opu(a)$ with $a\geq0$ ; a surjection
$\ft_{j-1}^*E\to\opu(z)$ (corresponding to the divisor $D_j$) and we
look for a surjection $\ft_{j-1}^*E\to\opu(y)$. We know that
$\ft_{j-1}(\pu)$ is not contained in $D_i$. There are a finite
number of points in $\pu$, say $x_1,\cdots,x_k$ such that
$\ft_{j-1}(x_l)\in D_i$. With our hypothesis on
$\a$ we have $y\geq z$ and $y\geq a-z$ (cf. proof of the preceding
proposition) so there always exists a section and we can choose it
such that the composition $\opu(a-z)\to f^*E\to\opu(y)$ is non zero
for $x_1,\cdots,x_k$. Then the new curve does not meet
$D_i\cap D_j$.

Because the condition is open we can find a curve for which it is true
for all $i$ and $j$.\hfill$\Box$

\begin{coro}
  Let $\a\in A_1(\Xt(\bar\phi))$ such that $\a\cdot \xi_k\geq0$ for
  all $k\in[1,n]$ and $\a\cdot\xi_x>0$ for some $x\in[1,n]$. 
Then there exists
  $\ft\in\Mor{\a}{\Xt(\bar\phi)}$ such that $\ft(\pu)$ meets $D_x$
  in $\displaystyle{D_x-\bigcup_{k\neq x}(D_x\cap D_k)}$.
\end{coro}

\dm
Let $\ft$ as in the preceding lemma. We know that $\ft(\pu)$ is not
contained in $D_x$ but has to meet $D_x$ (because of the
intersection number). The curve $\ft(\pu)$ does not meet any
intersection $D_i\cap D_j$, in particular it does not meet the
intersection $D_x\cap D_k$ for all $k$.\hfill$\Box$

\vs 0.4 cm

Let us now suppose that $\dim(X(\bar\phi))\geq3$. When
$\dim(X(\bar\phi))\leq2$ then $X(\bar\phi)$ is $\pu$ or $\p^2$ for
which 
the scheme of morphisms is well known. Let $D_x$ be a contracted
divisor, $\a\in A_1(\Xt(\bar\phi))$ and
$\ft\in\Mor{\a}{\Xt(\bar\phi)}$ as in the preceding corollary. There exists
$x_0\in\pu$ such that 
$$\ft(x_0)\in\displaystyle{D_x-\bigcup_{k\neq
    x}(D_x\cap D_k)}$$
and thanks to proposition \ref{recouvre}, for
any integer $i<x$ with $\sca{\a_i^\vee}{\a_x}=1$ there exists a curve 
$C\in[\Ga_{x,i}]$ such that $\ft(x_0)\in C$. 

\begin{prop}
\label{enunpoint}
Then there exists a deformation $\ft'$ of $\ft$ in
  $\Mor{\a}{\Xt(\bar\phi)}$ and an integer $i$ with
  $\sca{\a_i^\vee}{\a_x}=1$ such that $\ft'(\pu)$
  and $C$ meet exactly in $\ft(x_0)$ and transversaly.
\end{prop}

\dm
Let us first assume that $x<n$. 

\begin{lemm}
\label{existmoins1bis}
  There exists $j>x$ such that $\sca{\a_x^\vee}{\a_j}=1$.
\end{lemm}

\dm
It is enough to prove that there exists $j>x$ such that
$\sca{\b_x^\vee}{\b_j}\neq0$. Indeed taking the smallest such $j$ we must
have $\sca{\b_x^\vee}{\b_j}=-1$ because otherwise we would have
$\sca{\b_x^\vee}{\b_j}=2$ that is to say $\b_x=\b_j$. But for
$k\in[x+1,j-1]$ we have $\sca{\b_x^\vee}{\b_k}=0$ so in this case we have 
$$\bar\phi=s_{i(\b_1)}\cdots s_{i(\b_{x-1})}s_{i(\b_{x})}
s_{i(\b_{x+1})}\cdots s_{i(\b_{j-1})}s_{i(\b_{j})}s_{i(\b_{j+1})} \cdots
s_{i(\b_n)}$$
$$\ \ \ =s_{i(\b_1)}\cdots
s_{i(\b_{x-1})}s_{i(\b_{x})}s_{i(\b_{j})}
s_{i(\b_{x+1})}\cdots s_{i(\b_{j-1})}s_{i(\b_{j+1})} \cdots
s_{i(\b_n)}$$
$$= s_{i(\b_1)}\cdots
s_{i(\b_{x-1})}
s_{i(\b_{x+1})}\cdots s_{i(\b_{j-1})}s_{i(\b_{j+1})} \cdots
s_{i(\b_n)}\ \ \ \ \ \ \ \ \ \ $$ 
that is to say the expression $\bar\phi=s_{i(\b_1)}\cdots
s_{i(\b_{n})}$ was not reduced, a contradiction. Thus we have
$\sca{\b_x^\vee}{\b_j}=-1$. 
For such a $j$ we have $\sca{\b_x^\vee}{\b_k}=0$ for $k\in[x+1,j-1]$
and thus 
$$\sca{\a_x^\vee}{\a_j}=-\sca{\b_x^\vee}{\b_j}=1.$$

We have to prove that there exists $j>x$ such that
$\sca{\b_x^\vee}{\b_j}\neq0$. If not we would have:
$$\bar\phi=s_{i(\b_1)}\cdots s_{i(\b_{x-1})}s_{i(\b_{x})}
s_{i(\b_{x+1})}\cdots s_{i(\b_n)}= s_{i(\b_1)}\cdots
s_{i(\b_{x-1})}s_{i(\b_{x+1})}\cdots s_{i(\b_n)}s_{i(\b_{x})}$$
and we would have $\b_x=\b_n$ (remark \ref{debut}) thus
$\sca{\b_x^\vee}{\b_n}=2\neq0$ a contradiction.\hfill$\Box$

\vs 0.4 cm

In the case $x<n$ let $j$ be as in the lemma and consider the line
bundles $T_x$ and $T_j$. We have the formula (corollary \ref{tgt}):
$$T_i=\sum_{k=1}^{i}\sca{\a_k^\vee}{\a_i}\cdot\xi_k.$$
But $\sca{\a_k^\vee}{\a_i}\geq0$ for all $i$ and $k$
(proposition \ref{posi}) and $\a\cdot\xi_k\geq0$ for all $k$ by
assumption, therefore
$$\a\cdot T_x\geq\sca{\a_x^\vee}{\a_x}\a\cdot\xi_x= 2\a\cdot\xi_x>0\ \
{\rm and}\ \ \a\cdot
T_j\geq\sca{\a_x^\vee}{\a_j}\a\cdot\xi_x=\a\cdot\xi_x>0.$$

We construct the required $\ft'$ by induction on the fibrations. Let
us denote by $g:\pu\to\Xt(\bar\phi)$ the morphism whose image is $C$
(cf. proposition \ref{recouvre}) and define $P=\ft(x_0)$ and $P_k$ the
image of $P$ in $X_k/B_k$. Let us denote by $\ft_k$ (resp. $g_k$) the
morphism from $\pu$ to $X_k/B_k$ induced by $\ft$ (resp. by $g$). We
construct $\ft'$ by induction on the $\pu$-fibration beginning with
$\ft_{x-1}$.

\begin{lemm}
\label{diff}
  Let $\vp:X\to Y$ a $\pu$-bundle and $\a\in A_1(X)$ such that
  $\a\cdot T>0$ ($T$ is the relative tangent sheaf). Let
  $f\in\Mor{\a}{X}$ and $g:\pu\to X$ such that there exists
  $x_0\in\pu$ with $f(x_0)=g(x_0)$ and such that the images of $\vp\circ
  f$ and $\vp\circ g$ are distinct.

Then there exists a deformation $f'$ of $f$ meeting $g$ exactly in
$f(x_0)$ and transversaly.
\end{lemm}

\dm
Because the images are disctinct the curves $\vp\circ
  f(\pu)$ and $\vp\circ g(\pu)$ meet eachother in a finite number of
  points, say $x_0$ and $x_1,\cdots,x_k$. Let $E$ a rank 2 vector
  bundle defining the fibration, we can choose $E$ such that
  $(\vp\circ f)^*E=\opu\oplus\opu(a)$ with $a\geq0$. The morphism $f$
  is given by a surjection $s:(\vp\circ f)^*E\to\opu(y)$ with $\a\cdot
  T=2y-a>0$ (this implies $y>0$). A general surjection $s':(\vp\circ
  f)^*E\to\opu(y)$ will give a deformation of $f$. Because $y>0$,
  we can take such a surjection such that $s'(x_i)\neq s(x_i)$ for
  $i\in[1,k]$, $s'(x_0)=s(x_0)$ but are not equal at order 2 in
  $x_0$. This gives us a morphism $f'$ whose image meets the image of
  $g$ only in $f(x_0)$ and transversaly.\hfill$\Box$

\vs 0.4 cm

If $\ft_{x-1}(\pu)\neq g_{x-1}(\pu)$ then thanks to the lemma we can
construct $\ft'_{x}$ a deformation of $\ft_{x}$ meeting $g_{x}$
only in $P_x$. Taking by induction any section of $\ft'_{x}$ passing
trough the points $P_k$ for $k>x$ (this is possible because $\a\cdot
T_k\geq0$ for all $k$) we get the required deformation.

On the contrary if $\ft_{x-1}(\pu)=g_{x-1}(\pu)$ we use the following

\begin{lemm}
Let $\vp:X\to Y$ a $\pu$-bundle and $\a\in A_1(X)$ such that
  $\a\cdot T>0$ ($T$ is the relative tangent sheaf). Let
  $f\in\Mor{\a}{X}$ and $x_0\in\pu$.

There exists a deformation $f'$ of $f$ such that $f'$ and $f$
have distinct images still meeting in $f(x_0)$.
\end{lemm}

\dm
Let $E$ a rank 2 vector bundle defining the fibration, we can choose
$E$ such that $(\vp\circ f)^*E=\opu\oplus\opu(a)$ with $a\geq0$. The
morphism $f$ is given by a surjection $s:(\vp\circ f)^*E\to\opu(y)$
with $\a\cdot T=2y-a>0$ (this implies $y>0$). A general surjection
$s':(\vp\circ f)^*E\to\opu(y)$ will give a deformation of $f$. Because
$y>0$, we can take such a surjection such that $s'\neq s$ and
$s'(x_0)=s(x_0)$. This gives us the deformation $f'$.\hfill$\Box$

\vs 0.4 cm

If $\ft_x=g_x$ then, thanks to the lemma we can construct $\ft'_{x}$ a
deformation of $\ft_{x}$ meeting $g_{x}$ in $P_x$ and a finite number
of points. If $\ft_x\neq g_x$ we can take $\ft'_{x}=\ft_{x}$. Taking
by induction any section of $\ft'_{x}$ passing trough the points $P_k$
for $x<k<j$ (this is possible because $\a\cdot T_k\geq0$ for all $k$)
we get a deformation $\ft'_{j-1}$ of $\ft_{j-1}$ meeting $g_{j-1}$ in
$P_{j-1}$ and a finite number of points. Because we have $\a\cdot
T_j>0$ we can use lemma \ref{diff} to contruct a deformation $\ft'_j$
of $\ft_j$ meeting $g_j$ exactly in $P_j$ and transversaly. Taking by
induction any section of $\ft'_{j}$ passing trough the points $P_k$
for $k>j$ (this is possible because $\a\cdot T_k\geq0$ for all $k$) we
get the required deformation.

\vs 0.4 cm

The only case left is the case $x=n$. In this case, because $n\geq3$,
we can consider $\b_{n-1}$ and $\b_{n-2}$. Let us prove that
$\sca{\a_{n-2}^\vee}{\a_n}=\sca{\a_{n-1}^\vee}{\a_n}=1$. For
$\sca{\a_{n-1}^\vee}{\a_n}=-\sca{\b_{n-1}^\vee}{\b_n}=1$ it is just
corollary \ref{fin}. For $\sca{\a_{n-2}^\vee}{\a_n}$ we can apply
corollary \ref{fin} and it will be true except if
$\b_{n-2}=\b_n=\b$. But in this case we have
$$0\leq\sca{\a_{n-2}^\vee}{\a_n}= \sca{\b^\vee}{s_\b
  s_{\b_{n-1}}(\b_n)}=
\sca{\b^\vee}{s_\b(\b_{n-1}+\b_n)}= \sca{\b^\vee}{\b_{n-1}}=-1.$$
This is impossible and we must have
$\sca{\a_{n-2}^\vee}{\a_n}=\sca{\a_{n-1}^\vee}{\a_n}=1$. This in
particular implies that there are at least two $i<x=n$ such that
$\sca{\a_{i}^\vee}{\a_x}=1$, namely $i=n-1$ and $i=n-2$.

Let us now consider the morphism $\ft_{n-1}:\pu\to X_{n-1}/B_{n-1}$
induced by $\ft$ and two morphisms $g:\pu\to\Xt(\bar\phi)$ and
$h:\pu\to\Xt(\bar\phi)$ such that $g_*[\pu]=[\Ga_{n,n-1}]$ and
$h_*[\pu]=[\Ga_{n,n-2}]$. Because the classes $[\Ga_{n,n-1}]$ and
$[\Ga_{n,n-2}]$ are disctinct, the morphism $\ft_{n-1}$ has to be
distinct from one of the morphisms $g_{n-1}:\pu\to X_{n-1}/B_{n-1}$
and $h_{n-1}:\pu\to X_{n-1}/B_{n-1}$ deduced from $g$ and $h$. Let us
say that $\ft_{n-1}\neq g_{n-1}$ then applying lemma \ref{diff} we get
a deformation $\ft'$ of $\ft$ meeting $g$ only in
$\ft(x_0)$ and transversaly.\hfill$\Box$

\begin{prop}
\label{lissage}
  Let $C$ and $\ft'$ as in the preceding proposition, the curve
  $\ft'(\pu)\cup C$ can be smoothed. The smoothing is the image of a
  morphism $\fc:\pu\to\Xt(\bar\phi)$ and we have
$$[\fc(\pu)]\cdot\xi_x=[\ft'(\pu)]\cdot\xi_x-1=[\ft(\pu)]\cdot
\xi_x-1.$$
\end{prop}

\dm
We will use the following proposition proved in \cite{HH} corollary
1.2. for $\p^3$ but valid for any smooth variety $X$:

\begin{prop}
  Let $D$ be a nodal curve in a smooth variety $X$ and assume that
  $H^1(T_X\vert_D)=0$ then there exists a smooth deformation of $D$.
\end{prop}

In order to prove the proposition it suffices to prove that
$H^1(T_X\vert_D)=0$ where $X=\Xt(\bar\phi)$ and $D=\ft'(\pu)\cup C$
(which is a nodal curve). Let $P$ be the intersection point, we have
the exacte sequence
$$0\to\oo_{\ft'(\pu)}(-P)\to\oo_D\to\oo_C\to0$$
and it is enough to prove that $H^1(T_X\vert_C)=0$ and
$H^1(T_X\vert_{\ft'(\pu)}(-P))=0$.

One more time we do it by induction on the fibrations. Denote by
$\ft'_k:\pu\to X_k/B_k$ the morphism induced by $\ft'$ and
$C_k$ the image of $C$ in $X_k/B_k$. We assume that
$H^1(T_{X_{j-1}/B_{j-1}}\vert_{C_{j-1}})=0$ and
$H^1(T_{X_{j-1}/B_{j-1}}\vert_{\ft'_{j-1}(\pu)}(-P))=0$. We are going
to prove that $H^1(T_{X_{j}/B_{j}}\vert_{C_{j}})=0$ and
$H^1(T_{X_{j}/B_{j}}\vert_{\ft'_{j}(\pu)}(-P))=0$.

We have an exacte sequence 
$$0\to T_j\to T_{X_j/B_j}\to T_{X_{j-1}/B_{j-1}}\to 0$$
so it suffices to prove that $H^1(T_{j}\vert_{C_{j}})=0$ and
$H^1(T_{j}\vert_{\ft'_{j}(\pu)}(-P))=0$. But we have seen in the proof
of proposition \ref{recouvre} that $[\Ga_{x,i}]\cdot T_j\geq -1$ thus
the restriction of $T_j$ on $C_j$ is $\opu(u)$ with $u\geq-1$ and we
have the first vanishing. In the same way, we have $\a\cdot T_j\geq
0$ thus the restriction of $T_j$ on $\ft'_j(\pu)$ is $\opu(v)$ with
$v\geq0$ and we have the second vanishing.\hfill$\Box$

\section{Curves on minuscule Schubert varieties}

In this paragraph, we prove our main theorem on the irreducible
components of the scheme of morphisms from $\pu$ to $X(\bar\phi)$ a
minuscule Schubert variety.

\subsection{Moving out Schubert subvarieties}

We begin to prove that a general curve in $X(\bar\phi)$ is not
contained in a Schubert subvariety:

\begin{prop}
\label{movingout}
  Consider a morphism $f:\pu\to X(\bar\phi)$ such that $f$ factors
  through a Schubert variety $X(\bar\phi')\subset X(\bar\phi)$ (with
  $\bar\phi'<\bar\phi$) then there exists a deformation $f':\pu\to
  X(\bar\phi)$ of $f$ such that $f'$ does not factor trought
  $X(\bar\phi')$.
\end{prop}

\dm
Restricting ourselves to a smaller Schubert variety, we may assume
that the Schubert variety $X(\bar\phi')$ is a Schubert divisor $X(\bar
s_\b\phi)$ of $X(\bar\phi)$. But (theorem \ref{movingmin}) this
divisor has to be a moving divisor so $\b$ is simple and (proposition
\ref{movingdeb}) there exists a reduced expression
$$\bar\phi=s_{\ga_1}\cdots s_{\ga_n}$$
of $\bar\phi$ where $\ga_1=\b$ and $\bar\phi'=s_\b\bar\phi$. Consider
the Bott-Samelson resolution $\Xt(\bar\phi')$. If we denote with a
prime the corresponding elements in the Bott-Samelson construction we
have 
$B_{i+1}=s_{i(\b)}(B'_i)$, $P_{i+1}=s_{i(\b)}(P'_i)$ thus
$\Xt(\bar\phi')=X'_{n-1}/B'_{n-1}=s_{i(\b)}(f_n^{-1}\cdots
f_2^{-1}\sigma_1(X_0/B_0))$. 

This shows that we can identify $\Xt(\bar\phi')$ with the subscheme
$f_n^{-1}\cdots f_2^{-1}\sigma_1(X_0/B_0))$ of $\Xt(\bar\phi)$ which
is a fiber of the projection of $\Xt(\bar\phi)$ on $X_1/B_1$. We have
the commutative diagram:
$$\xymatrix{\Xt(\bar\phi')\ \ar[d]\ar@{^{(}->}[r]&\Xt(\bar\phi)\ar[d]\\
X(\bar\phi')\ \ar@{^{(}->}[r]&X(\bar\phi).}$$
The unipotent group $U_{-i(\b)}$ acts equivariantly on the second
vertical map and moves the first one. We may assume by induction that
$f$ does not factor through any Schubert subvariety of $X(\bar\phi')$
so we can find a section $g:\pu\to\Xt(\bar\phi')\subset\Xt(\bar\phi)$
of $f$. We can deform $g$ in $X(\bar\phi)$ thanks to the action of
$U_{-i(\b)}$ and we obtain a morphism  $g:\pu\to\Xt(\bar\phi)$ not
contained in $\Xt(\bar\phi')$. Projecting on $X(\bar\phi)$ gives a
deformation $f'$ of $f$ not contained in $X(\bar\phi')$.\hfill$\Box$ 

\begin{coro}
  For any morphism $f:\pu\to X(\bar\phi)$ there exists a deformation
  $f'$
of $f$ such that $f'$ does not factor trought any 
  $X(\bar\psi)\subset X(\bar\phi)$ (with $\bar\psi<\bar\phi$).
\label{sortir}
\end{coro}

\dm
Remark that if $f$ does not factor through a subvariety of
$X(\bar\phi)$ then it is also the case of any deformation. As there is
a finite number of Schubert varieties contained in $X(\bar\phi)$, we
apply the preceding proposition for each subvariety containing the
image of $f$.\hfill$\Box$

\vs 0.4 cm

Let $\pi:\Xt(\bar\phi)\to X(\bar\phi)$ a Bott-Samelson resolution. The
preceding result implies that for any morphism $f:\pu\to X(\bar\phi)$
there exist a deformation $f':\pu\to X(\bar\phi)$ of $f$ (the
deformation of the previous corollary) such that
$f'(\pu)$ meets the regular locus of $\pi$. We can thus consider a
section $\ft$ of $f'$:
$$\xymatrix{\pu\ar[dr]^{f'}\ar[r]^\ft&\Xt(\bar\phi)\ar[d]^\pi\\
&X(\bar\phi).}$$

\begin{rema}
The image $\ft(\pu)$ is not contained in any divisor $D_i$ on
$\Xt(\bar\phi)$. In fact, if it was the case it would means that
$f'(\pu)$ is contained in $\pi(D_i)$ which is a strict Schubert
subvariety of $X(\bar\phi)$. This is impossible.
\end{rema}

\begin{coro}
The morphism $\ft$ constructed from $f$ thanks to corollary
  \ref{sortir} is such that 
$$[\ft(\pu)]\cdot\xi_i\geq0\ \ {\rm for}\ \ {\rm all}\ \ i\in[1,n].$$
\end{coro}

\begin{prop}
  For any morphism $f:\pu\to X(\bar\phi)$ there exist a deformation
  $f'$ of $f$ such that $f'$ does not meet the image $\pi(D_x)$ of any
  contracted divisor $D_x$.
\end{prop}

\dm
We can replace $f$ by the deformation $f'$ of corollary
\ref{sortir}. We can thus assume that $f(\pu)$ is not contained in any
$\pi(D_i)$ for $i\in[1,n]$. We then have a section
$\ft:\pu\to\Xt(\bar\phi)$ of $f$. Let us denote $\a=\ft_*[(\pu)]$, we
have $\a\cdot\xi_i\geq0$ for all $i$. Define the subset
$A\subset\in[1,n]$ of all intergers $k$ such that $D_k$ is a
contracted divisor and set
$$l(\a)=\sum_{k\in A}\a\cdot\xi_k.$$
We prove the result by induction on $l(\a)$. If $l(\a)=0$ then $\ft$
does not meet any contracted divisor so $f$ does not meet the image
$\pi(D_x)$ of any contracted divisor $D_x$. 
Let $x$ be the smallest element in $A$ such that there exists a
morphism $g:\pu\to X(\bar\phi)$ not contained in any Schubert
subvariety with a section $\widetilde{g}:\pu\to\Xt(\bar\phi)$ such that
$\b=\widetilde{g}_*[(\pu)]$ with $l(\b)=l(\a)$, $\b\cdot\xi_x>0$ and
for which we have not constructed the required deformation yet. 
Thanks to propositions \ref{enunpoint}
and \ref{lissage}, for such a $g$ and $\widetilde{g}$ a section in
$\Xt(\bar\phi)$ there exists a deformation $\widetilde{g}'$ of
$\widetilde{g}$, an integer $i<x$ with $\sca{\a_i^\vee}{\a_x}=1$ and a
curve $C\in[\Ga_{x,i}]$ such that $\widetilde{g}'(\pu)\cup C$ can be
smoothed in $\widehat{g}(\pu)$. The morphism $\pi\circ\widehat{g}$
deforms to $g$ and we have
$\widehat{\b}=\widehat{g}_*[\pu]=\b+[\Ga_{x,i}]$. We thus have
$$\widehat{\b}\cdot\xi_k=
\left\{\begin{array}{cc}
\b\cdot\xi_k+1&{\rm for}\ k=i\\
\b\cdot\xi_k-1&{\rm for}\ k=x\\
\b\cdot\xi_k&{\rm otherwise}
\end{array}\right..$$
If $i\in A$ then because of our minimality assumption on $x$ we know
that there exists a required deformation for $\pi\circ\widehat{g}$ and
we can conclude because $g$ is a deformation of this deformation. If
on the contrary $i\not\in A$ then we have 
$$l(\widehat{\b})=\sum_{k\in A}\widehat{\b}\cdot\xi_k=l(\b)-1=l(\a)-1$$
and by induction there exists a required deformation for $\pi\circ\widehat{g}$
and we can conclude as above.\hfill$\Box$

\begin{coro}
\label{pointgen}
  There exists a dense open subset of $\Mor{}{X(\bar\phi)}$ whose
  elements do not meet the image $\pi(D_x)$ of any contracted divisor
  $D_x$ 
and are not contained in any $\pi(D_k)$ for $k\in[1,n]$.
\end{coro}

\begin{theo}
\label{final}
  Let $\a\in A_1(X(\bar\phi))$, the irreducible components of
  $\Mor{\a}{X(\bar\phi)}$ are indexed by $\comp(\a)$.
\end{theo}

\dm
We have a surjective morphism
$\pi_*:A_1(\Xt(\bar\phi))\to A_1(X(\bar\phi))$ and a natural morphism
$$\coprod_{\pi_*(\at)=\a}\Mor{\at}{\Xt(\bar\phi)}\to\Mor{\a}{X(\bar\phi)}.$$ 
Let us prove that the
irreducible components of $\Mor{\a}{X(\bar\phi)}$ are indexed by the
set $C(\a)$ of classes $\at\in A_1(\Xt(\bar\phi))$ such that
$\pi_*(\at)=\a$, $\at\cdot\xi_k\geq0$ for all $k\in[1,n]$ and
$\at\cdot\xi_x=0$ for all $x$ such that $D_x$ is a contracted divisor.

Because of corollary \ref{pointgen} we know that a general morphism
$f\in\Mor{\a}{X(\bar\phi)}$ can be lifted into
$\ft\in\Mor{\at}{\Xt(\bar\phi)}$ such that $\at\in C(\a)$. We thus
have a dominant morphism
$$\coprod_{\at\in C(\a)}\Mor{\at}{\Xt(\bar\phi)}\to\Mor{\a}{X(\bar\phi)}.$$
Let $\at\in C(\a)$ and $\ft$ a general element in
$\Mor{\at}{\Xt(\bar\phi)}$ (this scheme is irreducible thanks to
proposition \ref{irredenhaut}). We know (corollary \ref{pointgen} and
proposition \ref{irredenhaut}) that
its image is contained is the regular locus of $\pi$. If the morphism
$\pi\circ\ft$ was in the image of $\Mor{\at'}{\Xt(\bar\phi)}$ then
we would have a morphism $\ft'$ of class $\at'$ such that
$\pi\circ\ft'=\pi\circ\ft$. But because these curves are contained in the
regular locus of $\pi$ this implies that $\ft=\ft'$ and $\at'=\at$. The
images of the $\Mor{\at}{\Xt(\bar\phi)}$ for $\at\in C(\a)$ are the
irreducible components of $\Mor{\a}{X(\bar\phi)}$.

To conclude the proof we have to show that $C(\a)=\comp(\a)$. We begin
with the following

\begin{lemm}
  The kernel $K$ of the map $\pi_*:A_{n-1}(\Xt(\bar\phi))\to
  A_{n-1}(X(\bar\phi))$ is generated by the classes $\xi_x$ of the
  contracted divisors $D_x$.
\end{lemm}

\dm
M. Demazure proved in \cite{DE} that the morphism $\pi$ is an
isomorphism on the big cell of the Schubert variety
$X(\bar\phi)$. This in particular implies that the locus
$\widetilde{D}$ in $\Xt(\bar\phi)$ where $\pi$ is not an isomorphism is
contained in $\bigcup_iD_i.$
Moreover if the divisor $D_i$ is not contracted the open part
$D_i-\bigcup_{j\neq i}(D_i\cap D_j)$ is not contained in $\widetilde{D}$ so
that the codimension one part (in $\Xt(\bar\phi)$) of $\widetilde{D}$
is the union of the contracted divisors $D_x$.

Let us denote by $\widetilde{U}$ the open part in $\Xt(\bar\phi)$ where
$\pi$ is an isomorphism and $U$ its image in $X(\bar\phi)$. On the one
hand, the kernel of the surjective map $A_{n-1}(\Xt(\bar\phi))\to
A_{n-1}(\widetilde{U})$ is generated by the contracted divisors
$D_x$. On the other hand, we have
$A_{n-1}(\widetilde{U})=A_{n-1}({U})$ and because the complementary of
$U$ in $X(\bar\phi)$ is in codimension at least  2 (it is the image of
$\widetilde{D}$ with fibers of dimension at least 1 because Schubert
varieties are normal), we have
$A_{n-1}(U)=A_{n-1}(X(\bar\phi))$.\hfill$\Box$

\vs 0.4 cm

As $\Xt(\bar\phi))$ is smooth and projective we can identify
$A_{n-1}(\Xt(\bar\phi))^\vee$ with $A_{1}(\Xt(\bar\phi))$ and $\pi_*$
gives us a morphism
$A_{n-1}(\Xt(\bar\phi))^\vee\to A_{1}(X(\bar\phi))$. The
lemma leads to the following diagram whose first line is exact:
$$\xymatrix{0\ar[r]& A_{n-1}(X(\bar\phi))^\vee\ar[r]\ar[dr]^s&
  A_{n-1}(\Xt(\bar\phi))^\vee\ar[r]\ar[d]^\pi&K^\vee\ar[r]&0\\ 
&&A_{1}(X(\bar\phi))&&.}$$ 
Now we can translate the definition of $C(\a)$ in terms of
  $A_{n-1}(X(\bar\phi))^\vee$. Indeed, because of the wanishing
  condition on contracted divisor, all the elements of $C(\a)$ are in
  $A_{n-1}(X(\bar\phi))^\vee$ and go on $\a$ by $s$. What is left to
  prove is the following

  \begin{lemm}
An element $\at\in A_{n-1}(X(\bar\phi))^\vee$ seen as an element in
$A_{1}(\Xt(\bar\phi))$ is effective if and only if $\at\cdot\xi_i\geq0$ for
all $i\in[1,n]$. 
  \end{lemm}

\dm
We have seen proposition \ref{irredenhaut} that if all the
intersection $\at\cdot\xi_i$ are non negative then the class is
effective. 

Let $\at\in A_{n-1}(X(\bar\phi))^\vee$ an effective class. Because
  $\at$ is in $A_{n-1}(X(\bar\phi))^\vee$ we know that its
  intersection with all contracted $D_x$ are 0. Let $D_i$ a not
  contracted divisor, then its image in $X(\bar\phi)$ is a moving
  divisor (theorem \ref{movingmin}). Let $C$ a curve of class $\at$,
  if $C$ is not contained in $D_i$ then $C\cdot\xi_i\geq0$. If $C$
  is contained in $D_i$ then as in the proof of proposition
  \ref{movingout} we can deform this curve in the class $\at$ so
  that it is not contained in $D_i$ and we have
  $C\cdot\xi_i\geq0$.\hfill$\Box$

\vs 0.1 cm

This proves that $C(\a)=\comp(\a)$ and the theorem follows. Indeed,
$\comp(\a)$ is given (cf. paragraph \ref{preliminaires}) by the elements
$\b\in\pic(U)^\vee$ in the dual of the cone of effective divisors ($U$
is the dense orbit under ${\rm Stab}(X(\bar\phi))$). But
$\pic(U)=A_{n-1}(U)=A_{n-1}(X(\bar\phi))$ and the effective cone is
generated by the $\pi_*\xi_i$ with $D_i$ not contracted.\hfill$\Box$

\begin{small}

\vs 0.2 cm

\noi
{\textsc{Institut de Math{\'e}matiques de Jussieu}}

\vs -0.1 cm

\noi
{\textsc{175 rue du Chevaleret}}

\vs -0.1 cm

\noi
{\textsc{75013 Paris,}} \hs 0.2 cm{\textsc{France.}}

\vs -0.1 cm

\noi
{email: \texttt{nperrin@math.jussieu.fr}}

\end{small}


\begin{thebibliography}{FMcPSS}
\bibitem[Bo]{bourb} \textit{Nicolas Bourbaki}: {\'E}l{\'e}ments de
    math{\'e}matique. Fasc. XXXIV. Groupes et alg{\`e}bres de Lie. Chapitre
    IV : Groupes de Coxeter et syst{\`e}mes de Tits.
   Chapitre V: Groupes engendr{\'e}s par des r{\'e}flexions. Chapitre VI:
   syst{\`e}mes de racines. Actualit{\'e}s Scientifiques et
   Industrielles, No. 1337 Hermann, Paris 1968.
\bibitem[BP]{Brionpolo} \textit{Michel Brion and Patrick Polo}:
  Generic singularities of certain Schubert varieties. Math. Z. 231
  (1999), no. 2, 301--324.
\bibitem[De]{DE} \textit{Michel Demazure}: D{\'e}singularisation des
    vari{\'e}t{\'e}s de Schubert g{\'e}n{\'e}ralis{\'e}es. Collection of articles
    dedicated to Henri Cartan on the occasion of his 70th birthday, I.
    Ann. Sci. {\'E}cole Norm. Sup. (4) 7 (1974), 53--88.  
\bibitem[Fu]{F} \textit{William Fulton}: Intersection
  theory. Second edition. Ergebnisse der Mathematik und ihrer
  Grenzgebiete. 3. Folge. A Series of Modern Surveys in Mathematics,
  2. Springer-Verlag, Berlin (1998).
\bibitem[FMcPSS]{FS..} \textit{ William Fulton, Robert MacPherson, Frank
Sottile and Bernd Sturmfels}: Intersection theory on spherical
varieties. J. Algebraic Geom. 4 (1995), no. 1.
\bibitem[Ga]{Gaussent} \textit{Stephane Gaussent}: The fibre of the
  Bott-Samelson resolution. Indag. Math. (N.S.) 12 (2001), no. 4,
  453--468.
\bibitem[Gr]{GR}  \textit{Alexander Grothendieck}: Techniques de
  construction et th{\'e}or{\`e}mes d'existence en g{\'e}om{\'e}trie
  alg{\'e}brique. IV. Les sch{\'e}mas de Hilbert. (French). S{\'e}minaire
  Bourbaki, Vol. 6, Exp. No. 221, 249--276, Soc. Math. France, Paris,
  (1995).
\bibitem[HH]{HH} \textit{Robin Hartshorne and Andr{\'e} Hirschowitz}:
    Smoothing algebraic space curves. Algebraic geometry, Sitges
    (Barcelona), 98--131, L.N.M., 1124, Springer,
    Berlin, (1985).
\bibitem[K]{kempf} \textit{George R. Kempf}: Linear systems on
  homogeneous spaces. Ann. of Math. (2) 103 (1976), no. 3, 557--591.
\bibitem[LMS]{GofG/P3} \textit{Venkatramani Lakshmibai, Chitikila
  Musili and Conjeerveram S. Seshadri}: Geometry of
  $G/P$. III. Standard monomial theory for a quasi-minuscule
  $P$. Proc. Indian Acad. Sci. Sect. A Math. Sci. 88 (1979), no. 3,
  93--177.
\bibitem[LW]{LW} \textit{Venkatramani Lakshmibai and Jerzy Weyman}:
  Multiplicities of points on a Schubert variety in a minuscule
  $G/P$. Adv. Math. 84 (1990), no. 2, 179--208.
\bibitem[Ma]{maniv} \textit{Laurent Manivel}: Fonctions
  sym{\'e}triques, polyn{\^o}mes de Schubert et lieux de d{\'e}g{\'e}n{\'e}rescence. Cours
  Sp{\'e}cialis{\'e}s, 3. Soci{\'e}t{\'e} Math{\'e}matique de France, Paris, (1998).
\bibitem[Mo]{MO} \textit{Shigefumi Mori}: Projective manifolds with ample
    tangent bundles. Ann. of Math. (2) 110 (1979), no. 3, 593--606. 
\bibitem[P1]{PE} \textit{Nicolas Perrin}: Courbes rationnelles sur les
    vari{\'e}t{\'e}s homog{\`e}nes. Annales de l'Institut  Fourier, 52, no.1
    (2002), pp 105-132.
\bibitem[P2]{PE3} \textit{Nicolas Perrin}: Lieu singulier des surfaces
    rationnelles r{\'e}gl{\'e}es. Math. Z. 241 (2002),
   no. 2, 375--396. 
\bibitem[P3]{PE2} \textit{Nicolas Perrin}: Rational curves on
  homogeneous cones, in preparation. 
\end{thebibliography}
\end{document}